\documentclass[12pt,reqno]{amsart}
\usepackage[usenames]{color}
\usepackage{amsmath, amsfonts, epsfig, amssymb, graphics}
\usepackage{hyperref}
\def\charac{\raise 2pt\hbox{$\chi$}}
\def\Id{\mathrm{Id}}

\makeatletter
\def\revddots{\mathinner{\mkern1mu\raise\p@
\vbox{\kern7\p@\hbox{.}}\mkern2mu
\raise4\p@\hbox{.}\mkern2mu\raise7\p@\hbox{.}\mkern1mu}}
\makeatother

\makeatletter
\def\revots{\mathinner{\mkern1mu\raise\p@
\vbox{\kern7\p@\hbox{.}}\mkern2mu
\raise6\p@\hbox{.}\mkern2mu\raise12\p@\hbox{.}\mkern1mu}}
\makeatother

\def\Z{\mathbb{Z}}
\def\K{\mathbb{K}}
\def\L{\mathbb{L}}
\def\N{{\mathbb N}}

\def\A{{\mathcal A}}
\def\M{{M}}
\def\ox{\overline{x}\,}
\def\oy{\overline{y}\,}

\def\Ker{{\mathrm Ker}}

\def\carre(#1,#2)(#3){\put(#1,#2){\thicklines\framebox(#3,#3){}}}
\def\debut{{\rm s}}
\def\fin{{\rm e}}

\def\Paths{{\mathcal{P}}}
\def\NC{{\rm NoX}}

\def\bleu{\textcolor{blue}}

\def\rouge{\textcolor{red}}
\def\r#1{\case{#1}}
\def\dr#1{\dcase{#1}}
\def\etoile{\rouge{\star}}

\def\auteur#1{{\sc #1}}
\def\titreref#1{{\em #1}}
\def\vol#1{{\bf #1}}
\def\defn#1{{\bf #1}}

\def\jaune{\textcolor{yellow}}
\newdimen\carrelength
\def\jcarre{\jaune{\linethickness{\carrelength}\line(1,0){.85}}}
\def\case#1{\begin{picture}(1,1.2)(0,0)\setlength{\unitlength}{5mm}\put(-0.4,0.25){\jcarre}\put(-0.1,0){$$#1$$}\end{picture}}
\def\dcase#1{\begin{picture}(1,1.2)(0,0)\setlength{\unitlength}{6mm}\put(-0.4,0.25){\jcarre}\put(-0.3,0){$$#1$$}\end{picture}}
\def\Case#1#2{\put(#1){\case{#2}}}

\def\SL#1{{S\!L}_{#1}}
\def\monatop#1#2{\!\!\begin{array}{ll}\scriptstyle #1\\[-3pt] \scriptstyle #2\end{array}}

\def\N{{\mathbb N}}

\newtheorem{theorem}{Theorem}
\newtheorem{proposition}[theorem]{Proposition}

\newtheorem{corollary}[theorem]{Corollary}
\newtheorem{lemma}[theorem]{Lemma}

\def\X(#1,#2)#3#4{\multiput(#1,#2)(1,0){#4}{\O(0,0){#3}}}
\def\Y(#1,#2)#3#4{\multiput(#1,#2)(0,1){#4}{\O(0,0){#3}}}
\def\O(#1,#2)#3{\put(#1,#2)
                          {\put(0,0){\rouge{\circle*{.3}}}
                           \put(1,-1){\D(0,0){#3}}}}                         
\def\D(#1,#2)#3{\multiput(#1,#2)(1,-1){#3}{\bleu{\circle*{.3}}}}   
\def\pasX#1{\line(1,0){#1}}      
\def\pasY#1{\line(0,1){#1}} 
\def\coude(#1,#2)#3#4{\put(#1,#2){\pasX{#3}\pasY{#4}}}                          
\def\anticoude(#1,#2)#3#4{\put(#1,#2){\pasY{#3}\put(0,#3){\pasX{#4}}}}

\begin{document}

\title[Tilings]{\Large {\bleu{$\SL{k}$-Tilings of the Plane}}}
\author[F.~Bergeron and C.~Reutenauer]{Fran\c{c}ois Bergeron and Christophe Reutenauer}
\date{\today}

\address[F.~Bergeron and C.~Reutenauer]{D\'epartement de Math\'ematiques\\ Universit\'e
  du Qu\'ebec \`a Montr\'eal\\ 
 C.P. 8888\\ Succ. Centre-Ville\\ 
 Montr\'eal, Qu\'ebec, H3C 3P8, Canada.}
\email{bergeron.francois@uqam.ca, reutenauer.christophe@uqam.ca}

\thanks{F.~Bergeron and C.~Reutenauer are supported by NSERC-Canada  and FQRNT-Qu\'ebec.}

\maketitle
\begin{abstract}
We study properties of (bi-infinite) arrays having all adjacent $k\times k$ adjacent minors equal to one. 
If we further add the condition that all adjacent $(k-1)\times (k-1)$ minors be nonzero, then these arrays are necessarily of rank $k$. It follows that we can explicit construct all of them. Several nice properties are made apparent. In particular, we revisit, with this perspective, the notion of frieze patterns of Coxeter. This shed new light on their properties. A connexion is also established with the notion of $T$-systems of Statistical Physics.
\end{abstract}
 \parskip=0pt

{ \setcounter{tocdepth}{1}\parskip=0pt\footnotesize \tableofcontents}
\parskip=8pt

\section{Introduction}\label{introduction}
\begin{figure}[ht]
$$\begin{array}{ccccccccccccccc} \setlength{\unitlength}{3mm}
\ddots & \vdots & \vdots& \vdots& \vdots& \vdots& &\vdots & \vdots& \vdots& \vdots& \vdots& \vdots& \vdots &\revddots \\[3pt]
\cdots &887&567&247&174&101&\dr{28}&11&5&4&3&2&1&1&\cdots\\ 
\cdots &158&101&44&31&18&     \r{5}&2&1&1&1&1&1&2&\cdots \\
\cdots &61&39&17&12&7&           \r{2}&1&1&2&3&4&5&11&\cdots \\
\cdots &25&16&7&5&3&              \r{1}&1&2&5&8&11&14&31&\cdots \\
          &\dr{14}&\r{9}&\r{4}&\r{3}&\r{2}&\r{1}&\r{2}&\r{5}&\dr{13}&\dr{21}&\dr{29}&\dr{37}&\dr{82}&  \\
\cdots &3&2&1&1&1&\r{1}&3&8&21&34&47&60&133&\cdots \\
\cdots &1&1&1&2&3&\r{4}&13&35&92&149&206&263&583&\cdots \\
\cdots &1&2&3&7&11&\dr{15}&49&132&347&562&777&992&2199&\cdots \\
\revddots & \vdots & \vdots& \vdots& \vdots& \vdots& &\vdots    & \vdots& \vdots& \vdots& \vdots& \vdots& \vdots & \ddots
\end {array} $$
\begin{picture}(0,0)(0,0)\setlength{\unitlength}{3.5mm}
\put(-21,7.49){\rouge{\linethickness{1pt}\line(1,0){3}}}
\put(18,6.85){\rouge{\Huge$\longrightarrow$}}
\put(20,6.5){\rouge{$+$}}
\put(-2.6,14.5){\rouge{\linethickness{1pt}\line(0,1){2}}}
\put(-3.2,0.4){\rouge{\Huge$\downarrow$}}
\put(-2,0){\rouge{$+$}}
\end{picture}
\caption{A  $\SL{2}$-tiling with values in $\N^*$.}\label{fig1}
\end{figure}
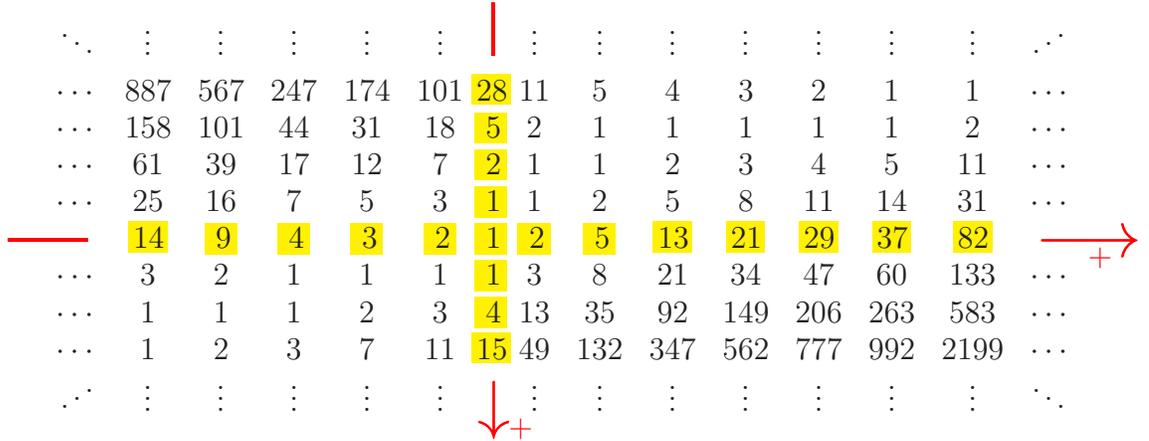

As discussed in \cite{assem}, the study of cluster algebras naturally leads to the special case $k=2$ of the notion of $\SL{k}$-tiling introduced in this paper. Our $\SL{k}$-tilings are simply $\Z\times \Z$ arrays of numbers (or elements of a commutative ring) having all adjacent $k\times k$ minors equal to one.  Not only are they a natural extension of notions already considered, but one can recast in their guise such notions as \defn{$T$-systems} of Theoretical Physics (see \cite{di_Francesco}), or \defn{frieze patterns} of Coxeter (see \cite{coxeter}).
An instance of a positive integer $\SL{2}$-tiling is given in Figure~\ref{fig1}.

Clearly any $\SL{k}$-tiling $\mathcal{A}$ is of rank at least $k$ (when considered as a bi-infinite matrix). As we will see, the $\SL{k}$-tilings that are of minimal rank are of particular interest, not only by themselves, but as well as for the cases when they correspond to frieze patterns or $T$-systems.  We call \defn{tame} such minimal rank $\SL{k}$-tilings, and we give several general results regarding them. Among these interesting results, we show that to any tame $\SL{k}$-tiling there corresponds another interesting tame $\SL{k}$-tiling, that we call its \defn{dual}. The entries of these dual tilings are obtained by computing adjacent $(k-1)\times (k-1)$-minors. It is striking that this duality is actually an involution.  We also re-derive, in a new an elegant manner, all the results of Conway-Coxeter concerning frieze-patterns. Indeed, our approach allows new tools to bear on this subject, especially because we can now make use of linear algebra and  particular presentations for $\SL{2}(\Z)$.

Our approach also opens the door for the study of generalized frieze patterns, including those that have already been considered in \cite{cordes}. This corresponds to the study of $\SL{k}$-tilings that afford two periods (in two linearly independent directions). We call \defn{toric} such $\SL{k}$-tilings, since they are evidently characterized by their value on a torus. Once again the tame situation  is of particular interest. More on this will be the subject of a planed sequel for this paper.
 
\section{Definitions}\label{sec_defn}
 We consider arrays $\A=(a_{ij})_{i,j\in\Z}$, with values $a_{ij}$ lying  
in a field $\K$. For  two equal cardinality finite subsets $I$ and $J$ of $\Z$, we denote by $\A_{IJ}$ the submatrix of $\A$ obtained by selecting the rows indexed by the elements of $I$ and columns indexed by the elements of $J$.  The corresponding minor is denoted\footnote{Without explicit reference to the underlying bi-infinite matrix $\mathcal{A}$.}  by $\M_{IJ}$, that is: $\M_{IJ}:=\det \A_{IJ}$. Since we often need to write down {\sl adjacent} $k\times k$ minors, we introduce the short hand notations:
\begin{equation}   \bleu{\A_{ij}^{(k)}:=\A_{\{i,\ldots, i+k-1\},\{j,\ldots, j+k-1\}}}, \qquad {\rm and}\qquad \bleu{\M_{ij}^{(k)}:=\det \A_{ij}^{(k)}}.\end{equation}

We say that $\A$ is a \defn{$\SL{k}$-tiling} of the plane if all its adjacent  $k\times k$ minors of $\A$ are equal to one. This is to say that it satisfies the \defn{$\SL{k}$-property}:
\begin{equation}\label{num_condition}
      \bleu{  \M_{ij}^{(k)}=1,\qquad  {\rm for\ all\ } i\  {\rm and\ } j\ {\rm in}\ \Z},
 \end{equation}
 We sometimes consider \defn{partial $\SL{k}$-tilings}, only defined on some subset  $S$ (called \defn{shape}) of $\Z\times \Z$, with condition~(\ref{num_condition}) applying only if all the entries considered belong to the underlying subset. As usual, a \defn{rectangle}  in $\Z\times \Z$ is a (possibly infinite) shape $S$ such that $(u,v+s)$ and $(u+r,v)$ lie in $S$, whenever $(u,v)$ and $(u+r,v+s)$ both lie in $S$. 
 A partial $\SL{k}$-tiling is said to be a \defn{$\SL{k}$-array} if its shape is a rectangle. In particular, $\N\times \N$ and $\Z\times \Z$ shaped $\SL{k}$-tilings are $\SL{k}$-arrays. Clearly linear combinations of rows (or columns) make sense for $\SL{k}$-arrays, so that we may consider the notion of rank of a such $\SL{k}$-tilings.  In particular any $\SL{k}$-array  is at least of rank $k$, since any $k$ consecutive rows have to be linearly independent in view of the $\SL{k}$-property. We say that a (partial) $\SL{k}$-tiling is  \defn{tame} if it has rank $k$. Otherwise we call it \defn{wild}.

A word of warning is in order concerning our convention for the underlying coordinate system. Indeed, as in the example of Figure~\ref{fig1}, we use the usual matrix convention for coordinates, so that the $x$-axis points downwards, and the $y$-axis points to the right.

\subsection*{A family of examples}
\begin{figure}[ht]
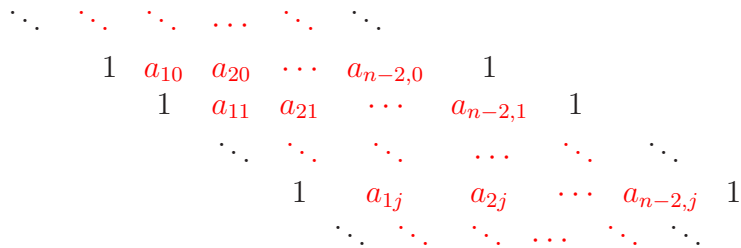

      $$\begin{array}{l}
         \begin{array}{ccccccccccccccccc}
                    &\, \ddots &\, \rouge{\ddots}&\, \rouge{\ddots} &\, \rouge{\cdots} &\, \rouge{\ddots} &\, \ddots&
        \end{array}\\[5pt]
         \begin{array}{ccccccccccccccccc}
         &\hbox{\hskip18pt}&      & 1& \rouge{a_{10}} & \rouge{a_{20}} &\rouge{\cdots} &  \rouge{a_{n-2,0}} & 1 &  \\
        &  &&       & 1& \rouge{a_{11}} & \rouge{a_{21}} &\rouge{\cdots} &  \rouge{a_{n-2,1}} & 1 & \\
         &&&&        &\,  \ddots &\, \rouge{\ddots} &\, \rouge{\ddots} &\, \rouge{\cdots} &\, \rouge{\ddots}&\,  \ddots &\\
          &&&&&      & 1& \rouge{a_{1j}} & \rouge{a_{2j}} &\rouge{\cdots} &  \rouge{a_{n-2,j}} & 1 & \\
       \end{array}\\[-2pt]
        \begin{array}{ccccccccccccccccc}
                  &\hbox{\hskip55pt}&&&&&&  & \ddots &  \rouge{\ddots}&\  \rouge{\ddots} &  \rouge{\cdots} &\  \rouge{\ddots} &  \ddots&
        \end{array}\\
       \end{array}$$ 
    \caption{Conway-Coxeter frieze patterns.}\label{fig1_5}
\end{figure}

The \defn{positive integer frieze patterns}  of Coxeter (see \cite{coxeter,conway,propp})  give rise to an interesting family of nonzero {\it partial} $\SL{2}$-tilings.  Up to a $45^\circ$ degree tilting, the original description of Coxeter may be formulated as follows. One considers partial $\SL{2}$-tilings\footnote{The $\SL{2}$ condition applies only when it makes sense.} such as the one illustrated in Figure~\ref{fig1_5}, assuming that all $a_{ij}$ are positive integers.
Note that the number of ``diagonals`` is $n$.

As shown in \cite{coxeter}, one of the striking property of frieze patterns is that they are necessarily periodic along the direction $y=x$. This is to say that there exists some $p$ in $\Z$ such that $a_{i+p,j+p}=a_{ij}$ for all $i$ and $j$, with $p=n+1$ ($n$ being the number of diagonals as above).

We may turn frieze patterns into full $\SL{2}$-tilings, by the simple device of extending them (skew) periodically both along the $x$ and $y$ directions, i.e.: setting 
    $$\bleu{a_{i+p,j}=-a_{ij}},\quad {\rm and}\quad \bleu{a_{i,j+p}=-a_{ij}}.$$
One needs only check that this is consistent with the $\SL{2}$-condition at the ``boundary''.
Such a tiling is illustrated  in Figure~\ref{fig2} in the case of a generic\footnote{All positive frieze patterns of width $2$ may be obtained from it by specialization.}   frieze pattern having $4$ diagonals. In this $\SL{2}$-tiling, $a$ and $b$ may assume any value as long as we have
     $${c=\frac{1+b}{a}},\qquad  {d=\frac{1+a+b}{a\,b}},\qquad {e=\frac{1+a}{b}}.$$
 Observe the further symmetry corresponding to a transposition followed by a diagonal translation.
\begin{figure}[ht]
      $$\begin {array}{rrrrrrrrrrrrrrrr} 
\rouge{a}&\rouge{b}&1&0&-1&-a&-b&-1&0&1&\rouge{a}&\rouge{b}\\ 
1&\rouge{c}&\rouge{d}&1&0&-1&-c&-d&-1&0&1&\rouge{c}\\ 
0&1&\rouge{e}&a&1&0&-1&-e&-a&-1&0&1\\ 
-1&0&1&b&c&1&0&-1&-b&-c&-1&0\\ 
-e&-1&0&1&d&e&1&0&-1&-d&-e&-1\\ 
-a&-b&-1&0&1&\rouge{a}&\rouge{b}&1&0&-1&-a&-b\\ 
-1&-c&-d&-1&0&1&\rouge{c}&\rouge{d}&1&0&-1&-c\\ 
0&-1&-e&-a&-1&0&1&\rouge{e}&a&1&0&-1\\ 
1&0&-1&-b&-c&-1&0&1&b&c&1&0\\ 
e&1&0&-1&-d&-e&-1&0&1&d&e&1\\ 
\rouge{a}&\rouge{b}&1&0&-1&-a&-b&-1&0&1&\rouge{a}&\rouge{b}\\ 
1&\rouge{c}&d&1&0&-1&-c&-d&-1&0&1&\rouge{c}\end {array}$$ 
    \caption{Skew-periodic extension of a  Conway-Coxeter frieze pattern.}\label{fig2}
\end{figure}
The number of  frieze patterns having $n$ diagonals, and for which all entries are positive integers, has been shown in \cite{conway} to be another incarnation of the ubiquitous Catalan numbers 
    $$C_{n-1}=\frac{1}{n}\binom{2\,n-2}{n-1}.$$
  A nice exposition of classical results regarding frieze patterns, as well as many new results tying their study to the type-$A$ cluster algebras of Fomin and Zelevinsky, is given by Propp in \cite{propp}.

It may readily be checked that the example of Figure~\ref{fig2} is a rank $2$ bi-infinite matrix. We will show in Section~\ref{sec_applications} how we may construct all frieze patterns using our theory of $\SL{k}$-tilings (for $k=2$), by extending them to complete $\SL{k}$-tilings. Moreover, using our theory, we give new proofs of all the results obtained by Coxeter and Conway.
Although our exploration of this point of view will mainly be for the case $k=2$, many of our results actually hold (with the necessary adaptations) in the general context of a suitable notion of $\SL{k}$-frieze patterns (see Section~\ref{sec_closing}).


\section{Tame $\SL{k}$-tilings}\label{sec_general}
Not all $\SL{k}$-tiling are tame, as seen in Example~(\ref{not_zero_free}) for $k=2$. 

\begin{equation}\label{not_zero_free}
 \begin {array}{cccccccccccr} 
\ddots&\vdots&\vdots&\vdots&\vdots&\vdots&\vdots&\vdots&\vdots&\vdots&\vdots&\revddots\\
\cdots&1&\bleu{x_{{11}}}&-1&\bleu{x_{{12}}}&1&\bleu{x_{{13}}}&-1&\bleu{x_{{14}}}&1&\bleu{x_{{15}}}&\cdots\\ 
\cdots&0&1&0&-1&0&1&0&-1&0&1&\cdots\\ 
\cdots&-1&\bleu{x_{{21}}}&1&\bleu{x_{{22}}}&-1&\bleu{x_{{23}}}&1&\bleu{x_{{24}}}&-1&\bleu{x_{{25}}}&\cdots\\ 
\cdots&0&-1&0&1&0&-1&0&1&0&-1&\cdots\\ 
\cdots&1&\bleu{x_{{31}}}&-1&\bleu{x_{{32}}}&1&\bleu{x_{{33}}}&-1&\bleu{x_{{34}}}&1&\bleu{x_{{35}}}&\cdots\\ 
\cdots&0&1&0&-1&0&1&0&-1&0&1&\cdots\\ 
\cdots&-1&\bleu{x_{{41}}}&1&\bleu{x_{{42}}}&-1&\bleu{x_{{43}}}&1&\bleu{x_{{44}}}&-1&\bleu{x_{{45}}}&\cdots\\ 
\cdots&0&-1&0&1&0&-1&0&1&0&-1&\cdots\\ 
\cdots&1&\bleu{x_{{51}}}&-1&\bleu{x_{{52}}}&1&\bleu{x_{{53}}}&-1&\bleu{x_{{54}}}&1&\bleu{x_{{55}}}&\cdots\\ 
\cdots&0&1&0&-1&0&1&0&-1&0&1&\cdots\\
\revddots&\vdots&\vdots&\vdots&\vdots&\vdots&\vdots&\vdots&\vdots&\vdots&\vdots&\ddots\end {array}
\end{equation}
Here the $x_{ij}$ may be chosen at will (or as independent variables). In particular, this example shows that there are $\SL{2}$-tilings of any rank $\geq 2$.

In part, the interest of considering tame tilings comes from the fact that they are easily characterized by their value on relatively small subsets of $\Z\times\Z$. But we will also make evident that tame tilings have very nice properties.
We first illustrate tameness with the following special case.

\subsection*{$0$-free tilings} We say that a $\SL{k}$-tiling is  \defn{$0$-free} if all its $(k-1)\times (k-1)$ adjacent subminors are nonzero. Note that in the case $k=2$: a $\SL{2}$-tiling is 0-free if its values are nonzero; in particular if they are positive integers, as in Figure~\ref{fig1} or the $\SL{2}$-tilings constructed in \cite{assem}.

The proof of the following proposition uses Dodgson\footnote{a.k.a. Lewis Carrol} ``Condensation Law of Determinants'' \cite{alice}, that can be stated in the format:
\begin{equation}\label{carrol}
   \bleu{\M_{ij}^{(r+1)} \M_{i+1,j+1}^{(r-1)} = \det
         \left(\begin{array}{lll} \M_{ij}^{(r)}  & \M_{i,j+1}^{(r)}\\[6pt]
                                 \M_{i+1,j}^{(r)}  &    \M_{i+1,j+1}^{(r)}
          \end{array} \right)}
  \end{equation}
for all $r$.
In fact this is a direct consequence of a result of Desnanot and Jacobi (see \cite[Th. 3.12, page 111]{bressoud}).
For instance, with $r=2$, we get the identity
$$\begin{array}{l}
    \det 
         \left(\begin{array} {lll}a_{ij}  & a_{i,j+1} & a_{i,j+2}\\
                                  a_{i+1,j}  & a_{i+1,j+1} & a_{i+1,j+2}\\
                                  a_{i+2,j}  & a_{i+2,j+1} & a_{i+2,j+2}\\
              \end{array}\right) \det 
         \begin{pmatrix} a_{i+1,j+1}
          \end{pmatrix} =\\[25pt] 
          \hskip150pt \det \begin{pmatrix} \left| \begin{array}{llll} a_{ij}    & a_{i,j+1}\\  
                                                                         a_{i+1,j}&a_{i+1,j+1}
                                                   \end{array}\right| &
                                           \left| \begin{array}{llll}  a_{i,j+1}    & a_{i,j+2}\\ 
                                                                         a_{i+1,j+1} & a_{i+1,j+2}
                                                   \end{array}\right| \\[14pt]
                                           \left| \begin{array}{llll}  a_{i+1,j}    & a_{i+1,j+1}\\ 
                                                                         a_{i+2,j}& a_{i+2,j+1}
                                                  \end{array}\right| &
                                           \left| \begin{array}{llll}  a_{i+1,j+1}& a_{i+1,j+2}\\ 
                                                                              a_{i+2,j+1}& a_{i+2,j+2}
                                                 \end{array}\right| 
          \end{pmatrix}
          \end{array}$$
          
     \begin{proposition}\label{rankk}
   \bleu{Any $0$-free $\SL{k}$-array is tame}.
\end{proposition}

 \begin{proof}[\bf Proof.]
Consider any adjacent $(k-1)\times (k-1)$ subarray of a  $\SL{k}$-array, we observe that the determinant of the corresponding submatrix does not vanish, since this is precisely the $0$-free condition. On the other hand, for $r=k$,  the right-hand side of (\ref{carrol}) is zero in all instances, since the four $k\times k$ minors considered are all equal to $1$. We thus conclude that any adjacent $(k+1)\times (k+1)$ subarray of a $0$-free $\SL{k}$-array must necessarily have vanishing determinant. The proof then follows from Lemma~\ref{det_rank} (see Section~\ref{preuves}).
  \end{proof}

We may construct all tame $\SL{k}$-tilings as follows. Given a $\SL{k}$-tiling $\A$, let us denote by $R_i$ and $C_j$ its rows and columns. Then each $C_j$ is a linear combination of the $k$ preceding columns, that is, of $C_{j-1},\ldots, C_{j-k}$. The linear combination may be written as 
\begin{equation}\label{colonnes}
   \bleu{ (-1)^kC_{0}-(-1)^{k}a_1C_{1}+\cdots -a_{k-1}C_{k-1}+C_k=0}.
 \end{equation}
Indeed this follows from the $\SL{k}$-property and from the next lemma, which is an exercise in linear algebra (expansion of the determinant with respect to the rows), left to the reader.

\begin{lemma}
\bleu{Let $A$ be a rank $k$ matrix with $k+1$ columns $C_0,\ldots,C_k$ and rows indexed by $\Z$. Then}
\begin{equation}\begin{array}{ll}
   \bleu{(-1)^k M_{01}^{(k)}C_{0}-(-1)^{k}\det A_{\{0,\ldots,k-1\},\{0,2,\ldots,k-1\}} C_{1}+\cdots}\\[6pt]
\qquad \qquad\qquad\qquad\qquad
 \bleu{-\det A_{\{0,\ldots,k-1\},\{0,\ldots,k-2,k\}} C_{k-1}+ M_{00}^{(k)}C_k=0}.
\end{array}
\end{equation}
\end{lemma}

Hence, to each $j$ in $\Z$ we associate a row vector $\lambda_j=(a_1,...,a_{k-1})$ of dimension $m=k-1$ over the field $\K$, and we simply denote by $\lambda$ the resulting element of $ (\K^{1\times m})^{\Z}$. Similarly, we associate to each $i$ a column vector $\gamma_i$ of dimension $m$ over $\K$, which expresses the linear dependence of $R_i$ on the $m$ preceding rows, and denote by $\gamma$ the resulting element of $(\gamma_i)\in (\K^{m\times 1})^{\Z}$. These row and columns vectors are called the \defn{linearization coefficients} of the $\SL{k}$-tiling $\A$. 
We call \defn{linearization data} the triple 
\begin{equation}\label{data_lin}
   \bleu{(\A_{00}^{(k)}, \lambda, \gamma)} \quad {\rm in}\quad \bleu{ \SL{k}(\K)\times  (\K^{1\times m})^{\Z} \times (\K^{m\times 1})^{\Z}},\quad \rm {where}\ \bleu{m=k-1}.
 \end{equation}

\begin{proposition}\label{coefflinear}
\bleu{The mapping}
     $$\bleu{\A \mapsto (\A_{00}^{(k)}, \lambda, \gamma)},$$
\bleu{which associates to a tame $\SL{k}$-tiling  its linearization data, is a bijection between the set of tame $\SL{k}$-tilings and the set $\SL{k}(\K)\times  (\K^{1\times m})^{\Z} \times (\K^{m\times 1})^{\Z}$, with $m=k-1$}.
\end{proposition}

\begin{proof}[\bf Proof.]
The fact that this mapping is well-defined and injective follows from the remarks preceding the proof. For surjectivity, let the data in $\SL{k}(\K)\times  (\K^{1\times m})^{\Z} \times (\K^{m\times 1})^{\Z}$ be given. Then clearly there exists a $\Z \times \Z$ array $\A$ of rank $k$ which maps onto this data. 
We have only to verify that the $\SL{k}$-property holds. This is a consequence of the following easy fact (and its dual): let $M$ be a $k\times (k+1)$-matrix with columns $C_0,\ldots,C_k$ such that (\ref{colonnes}) holds.
Then the matrix of its first $k$ columns is of determinant 1 if and only if the matrix of its $k$ last columns is of determinant 1.
Indeed, 
\begin{eqnarray*}
\det(C_1,\ldots,C_{k})&=&\det(C_1,\ldots ,C_{k-1},-(-1)^kC_{0}+(-1)^{k}a_1C_{1}-\ldots +a_{k-1}C_{k-1})\\
  &=&\det(C_1,\ldots ,C_{k-1},-(-1)^kC_{0})\\
  &=&\det(C_{0},\ldots,C_{k-1}).
 \end{eqnarray*}
\end{proof}

\begin{proposition}\label{rank1}
\bleu{Let $\A$ be an $\SL{k}$-tiling. Then $\A$ is tame if and only if the infinite matrix $\mathcal{M}:=(\M_{IJ})_{I,J}$,
with $I$ and $J$ varying in $k$-subsets of $\Z$, is of rank $1$. In particular, if $I_0, J_0$ are intervals in $\Z$ and $I, J$ are any $k$-subsets, then we have} 
\begin{equation}\label{mineursmultiplicatifs}
    \bleu{\M_{IJ}=\M_{IJ_0} \, \M_{I_0J}}.
  \end{equation}
\end{proposition}

\begin{proof}[\bf Proof.]
If $\A$ is tame, then the fact that $\mathcal{M}$ has rank $1$ is a consequence of the study of $\N\times\N$ tilings in Section~\ref{NbyN}. The converse follows from Lemma~\ref{rang1}. 

Now take $I$, $J$, $I_0$, and $J_0$ as in the statement. Then
\begin{equation}\label{etape}
\bleu{
\begin{pmatrix}
\M_{I_0J_0} & \M_{I_0J} \\
\M_{IJ_0} & \M_{IJ}
\end{pmatrix}
}
\end{equation}
is a submatrix of $\mathcal{M}$. But $\M_{I_0J_0} =1$ (since $\A$ is a $\SL{k}$-tiling), and the determinant of (\ref{etape})  mush vanish, i.e.:
      $$\M_{IJ}-\M_{I_0J}\M_{IJ_0} =0,$$
since $\mathcal{M}$ is of rank $1$.
Thus we get the desired equality.
\end{proof}

The easy direct proof of the next result is left to the reader. It will be useful in the sequel.

\begin{lemma}\label{suppression_colonne}
\bleu{Let $C_1$, $C_2$, and $C_3$ be three consecutive columns {\rm (}resp. rows{\rm )}, of a tame $\SL{2}$-tiling, that are such that $C_2=C_1+C_3$. Then a new tame $\SL{2}$-tiling may be constructed by suppressing the column {\rm (}resp. row{\rm )}  $C_2$}.
\end{lemma}

\subsection*{Group actions on tilings}
There is a natural translation action of $\Z^2$ on $\SL{2}$-tilings.
Formally, the action of the vector $(p,q)$ replaces the tiling $\A=(a_{ij})$ by 
    $$\bleu{(p,q)\cdot \A=(a_{i+p,j+q})}.$$
 We denote by $\A_x$ the translate of $\A$ by $(1,0)$, and by $\A_y$ the translate by $(0,1)$.
We may describe these last translates intrinsically in terms of the data given by the bijection of Proposition~\ref{coefflinear}.
More specifically, let $(\mathcal{S}, \lambda, \gamma)$  be the linearization data corresponding to $\A$ via this bijection. Then the linearization data corresponding to  $\A_x$ is  $(\mathcal{S}_x, \lambda, \gamma')$, with
$$
\mathcal{S}_x = \begin{pmatrix} 0&1&0 &\ldots&0 \\ 0&0&1&\ldots&0 \\ \vdots&\ddots&\ddots &\ddots&\vdots\\ 0&\cdots&0&0&1\\(-1)^{k-1} & (-1)^k(\gamma_k)_1&\ldots&&(\gamma_k)_{k-1}
\end{pmatrix} \mathcal{S},
$$
and $\gamma_i'=\gamma_{i+1}$
Transposition of matrices (which amounts to exchanging rows and columns) also preserves $\SL{k}$-tilings, as is easily observed. In terms of linearization data, it amounts to transposing the initial matrix and exchanging $\lambda$ and $\gamma$. We can thus easily describe $\A_y$ using these remarks.

Observe that if $k\equiv 0$ or $1\mod 4$, then a vertical or an horizontal symmetry also preserves $\SL{k}$-tilings, since in this case, such a symmetry preserves the determinant of $k\times k$ matrices (for $k\equiv 2$ or $3 \mod 4$, the determinant is clearly multiplied by $-1$).

\section{Dual tilings} 
To any array $\A$, we associate the \defn{$m$-derived array}:
 \begin{equation}\label{rule}
    \bleu{\partial_m \A:=\left(\M_{ij}^{(m)}\right)_{i,j}},
\end{equation}
consisting of the adjacent $m\times m$ minors of $\A$. For $\SL{k}$-arrays, we are specially interested in the case  $m=k-1$, in which case the resulting array is called the \defn{dual} array of $\A$. We also write $\A^*$ for $\partial_{k-1} \A$.
Clearly $\partial_ 1$ is the identity operator, and it is natural to set $\partial_0 \A$ equal to the tiling whose value is $1$ in all positions.
As an illustration of the above definition, the dual of the $\SL{3}$-tiling 
\begin{equation}\label{pavage_fibo3}
\bleu{\begin {array}{ccccccccccccccc} 
\ddots&\vdots & \vdots &\vdots &\vdots &\vdots &\vdots &\vdots &\vdots &\revddots\\
\cdots & 8997&1782&353&70&14&3&1&1&\cdots \\ 
\cdots & 1782&353&70&14&3&1&1&2&\cdots\\ 
\cdots & 353&70&14&3&1&1&2&5&\cdots\\ 
\cdots & 70&14&3&1&1&2&5&14&\cdots\\ 
\cdots & 14&3&1&1&2&5&14&42&\cdots\\ 
\cdots & 3&1&1&2&5&14&42&131&\cdots\\ 
\cdots & 1&1&2&5&14&42&131&417&\cdots\\
\revddots&\vdots & \vdots &\vdots &\vdots &\vdots &\vdots &\vdots &\vdots &\ddots
\end {array}}
\end{equation}
 is the tiling
\begin{equation}\label{une_derivee}
\bleu{ \begin {array}{ccccccccccccccccc} 
\ddots&\vdots & \vdots &\vdots &\vdots &\vdots &\vdots &\vdots &\revddots\\
\cdots & 417&131&42&14&5&2&1&\cdots\\ 
 \cdots &131&42&14&5&2&1&1&\cdots\\ 
\cdots & 42&14&5&2&1&1&3&\cdots\\ 
\cdots &14&5&2&1&1&3&14&\cdots\\ 
\cdots &5&2&1&1&3&14&70&\cdots\\ 
\cdots &2&1&1&3&14&70&353&\cdots\\
\revddots&\vdots & \vdots &\vdots &\vdots &\vdots &\vdots &\vdots &\ddots
\end {array} }
\end{equation}
 





We have the following property of derivation of arrays, that will be proved in Section~\ref{preuves}.

\begin{proposition}\label{derive_equation}
\bleu{The dual of a tame $\SL{k}$-tiling is a tame $\SL{k}$-tiling}.
  \bleu{Moreover, for any natural integers $r,s$ such that $r+s=k$, we have}
 \begin{equation}\label{condensation_devive}
   \bleu{\partial_r \A^*=(r-1,r-1)\cdot (\partial_{s} \A)}
 \end{equation}
\bleu{In particular, $(\A^*)^*$ and $\A$  coincide up to translation}.
\end{proposition}

Observe also that, with $r=k-1$,  identity (\ref{carrol}) gives $(\A^*)^*=\rouge{a_{11}}\det(\A)$ for any $3\times 3$ matrix
 $${\A=\begin{pmatrix} 
              a_{00}& a_{01} & a_{02}\\
              a_{10}& \rouge{a_{11}} & a_{12}\\
              a_{20}& a_{21} & a_{22}
           \end{pmatrix}}.$$
Hence, for any $\SL{3}$-tiling $\A$, we have   
 \begin{equation}\label{rec_trois}
     \bleu{ (\A^*)^* = \left(a_{i+1,j+1}\, \M_{ij}^{(3)}\right)_{i,j\in \Z}}.
   \end{equation}
It follows that $(\A^*)^*=\A$ (up to the necessary translation) for any  $\SL{3}$-tiling, wether they be tame or wild. However, for $k\geq 4$, it may be checked that the tameness condition is necessary for (\ref{condensation_devive}) to hold.  

\section{Tilings associated to paths}\label{words_tilings}
We consider  \defn{paths} $\pi$
as lists of \defn{points}, i.e.: elements of $\Z\times \Z$,
    $$\pi=(i_0,j_0),(i_1,j_1),\ldots, (i_{N},j_{N}),$$
starting at $\debut(\pi):=(i_0,j_0)$ and ending at $\fin(\pi):=(i_{N},j_{N})$, and such that
    $$(i_{s+1},j_{s+1})=\begin{cases} (i_s,j_s-(1,0) & \text{or}, \\[6pt]
             (i_s,j_s)+(0,1),
\end{cases}$$
for points along the path. To understand why signs appear here, it may be good to recall our convention for the orientation of the $x$ and $y$ axis (see Section~\ref{sec_defn}).
If we fix the start and end points $(i_0,j_0)$ and $(i_{N},j_{N})$, it is well known that these paths number $\binom{m+n}{m}$ (with $(-m,n)=(i_{N},j_{N})-(i_0,j_0)$), and that they are in bijection with words  
\begin{equation}\label{le_mot}
    w=w_1w_2\cdots w_{n+m}
\end{equation}
on the alphabet $\A=\{x,y\}$,  having $m$ occurrences of the letter $x$, and $n$ occurrences of the letter $y$. Recall that the corresponding classical bijection, between paths and words, is realized by choosing $(i_{s+1},j_{s+1})=(i_s,j_s)-(1,0)$ if $w_s=x$, and $(i_{s+1},j_{s+1})=(i_s,j_s)+(0,1)$ if $w_s=y$. We denote by $\pi_w$ the resulting path. 

We also consider words, and associated paths, that are infinite in both directions, 
 \begin{equation}\label{le_mot}
    w=\cdots w_{-3}w_{-2}w_{-1} w_0w_1w_2w_3 \cdots
\end{equation}  
and say that they are \defn{bi-infinite} words (or paths).  Such a word (and the associated path) is said to be \defn {admissible} if there are infinitely many $x$' s and $y$'s in both directions.

Let us now associate to a given word $w$ (finite or bi-infinite) a tiling $\A=\A_{w;k}$, whose entries are obtained by the (weighted) enumeration of paths starting and ending at some points of $\pi_w$. We restrict these paths to stay within some ``distance'' $k$ of $\pi_w$. This is made more precisely below after the introduction of more notation and terminology (some of which will only be used later).

\begin{figure}[ht]\setlength{\unitlength}{5mm}
 \begin{picture}(0,15)(-1.5,-8)
 \put(-17,0){$\rouge{x}$}
 \put(-12,5.5){$\rouge{y}$} 
 \put(-16,5){\rouge{\vector(0,-1){10}}}
  \put(-16,5){\rouge{\vector(1,0){10}}}
 \put(6,-6.5){\linethickness{4pt}\jaune{\line(0,1){10.5}}}
  \put(5.5,-7){\linethickness{4pt}\jaune{\line(-1,0){13.5}}}
 \put(-7,-7){\coude(-1,0)41\coude(3,1)33\coude(6,4)22\coude(8,6)11\coude(9,7)44}
   \X(-10,-8){0}{2}
   \Y(-8,-8){0}{1}
   \Y(-8,-7){1}{1}
   \Y(-7,-7){1}{1}
   \Y(-7,-6){2}{1}
   \X(-7,-5){2}{4}
   \Y(-3,-5){2}{2}
   \X(-3,-3){2}{3}
   \Y(0,-3){2}{3}
   \X(0,0){2}{4}
   \Y(4,0){2}{1}
    \X(4,1){2}{3}
    \Y(6,2){2}{3}
    \X(6,4){2}{3}
    \X(9,4){1}{1}
    \X(10,4){0}{1}
    \put(5.8,-7){$p=(i,j)$}
      \put(-16,-7.2){$(i,\beta_\pi(p))=\gamma_w(p)\rightarrow$}
      \put(5,5.5){$\chi_w(p)=(\alpha_\pi(p),j)$}
      \put(5.8,4.5){$\downarrow$}
          \put(12.5,4.5){\rouge{$\pi_w$}}
                    \put(10.8,4.2){\rouge{\vector(3,1){1}}}
          \put(-13,-9){\rouge{$\pi_w$}}
                    \put(-11.8,-8.4){\rouge{\vector(3,1){1}}}
            \end{picture}
\caption{Path from $\gamma_w(p)$ to $\chi_w(p)$ in the $3$-fringe of $\pi_w$.}\label{fig_dessous}
\end{figure}
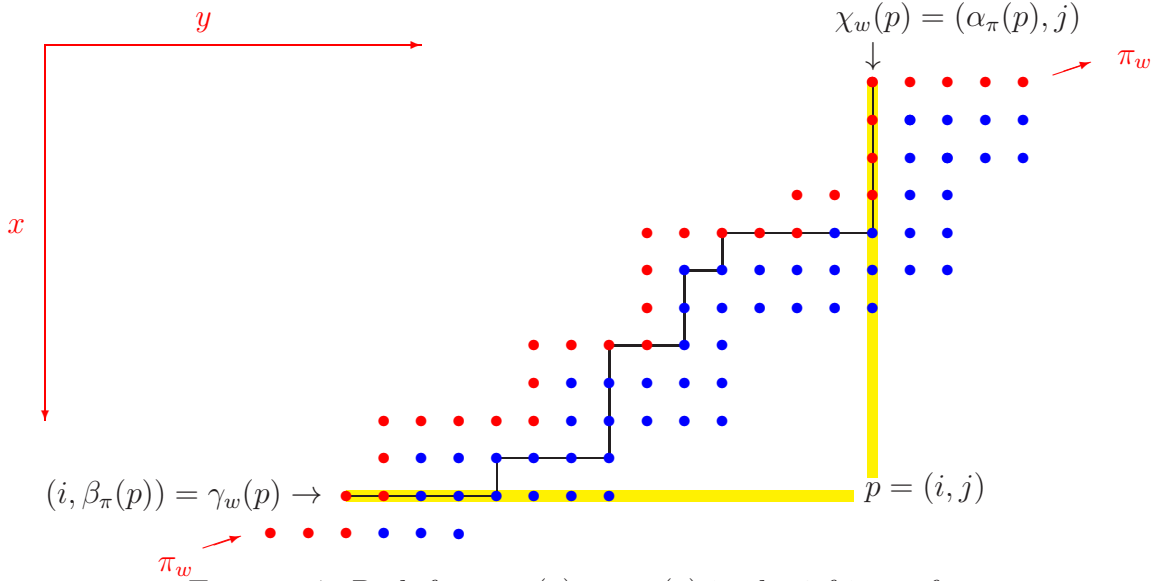 

Given a path $\pi_w$ as above, for each point $p=(i,j)$  in $\Z\times \Z$ we denote by $\gamma_w(p)=(i,\beta_\pi(p))$ (resp. $\chi_w(p)=(\alpha_\pi(p),j)$) the leftmost (resp. topmost) point 
that lies on the path $\pi_w$, which has the same first (resp. second) coordinate as $p$. We say that these are respectively the \defn{horizontal projection} and \defn{vertical projection} of $p$ on $\pi_w$. 
A point $p=(i,j)$ is said to lie \defn{below} the path $\pi$ if we have the inequalities $\beta_\pi(p)\leq j$, or equivalently, $\alpha_\pi(p)\leq i$.  Otherwise, we say that $p$ lies \defn{above} the path.  

We now consider the word $w(p)$:
associated to the portion of the path $\pi_w$ going from the horizontal projection $\gamma_w(p)$ of $p$ on $\pi_w$, to its vertical projection $\chi_w(p)$. This word is used to define the notion of \defn{projection word} of a point $p$, denoted by $w_p$, as follows. We simply set $w_p:=w(p)$ whenever $p$ lies below the path. Otherwise, when $p$ lies above the path, we set $w_p:=\overline{w(p)}$. Here $\overline{w}$ is the operation corresponding to reading the letters of a word $w$ in reverse order, replacing each $x$ by $\overline{x}$ and each $y$ by $\overline{y}$. For example, with $w=yyxyxyyyx$, we get $\overline{w}=\ox\oy\oy\ox\oy\ox\ox\oy\oy$.

\begin{figure}[ht]
 \begin{center}
\begin{picture}(-50,130)(0,0)
 \put(-160,70){\rouge{$\begin {array}{ccccccccc} 
 \ddots&\vdots&\vdots&\vdots&\vdots&\ \vdots\ &\ \, \vdots\ \,&\quad \vdots\quad &\revddots\\
\cdots&  1437&457&30&10\\ 
\cdots&   457&146&10&\ \  4\ \  \\ 
\cdots&   30&10\\ 
\cdots&   10&\ 4\ \\ 
\cdots\\ 
\cdots\\ 
\cdots
\end {array}$}}
\put(-150,55){\bleu{$\begin {array}{cccccccccc}  
 &&&&&&&& \\
 && & & & \r{1}&\ \r{1}\ &\ \r{1}\ &\cdots\\  
& & & & & \r{1}& 2&3&\cdots\\ 
& & &\r{1}&\r{1}&\r{1}&3&\ 6&\cdots\\  
& & &\r{1}&2&3& 10&22&\cdots\\ 
&\ \r{1}\ &\r{1}&\r{1}&3&6&\ 22\  &53&\cdots\\  
&\r{1}& 2&3&10&22&84&211&\cdots\\ 
& \r{1}&\quad 3\quad &6&\ 22\ & 53&  211 &553&\cdots\\
\revddots&\vdots&\vdots&\vdots&\vdots&\vdots&\vdots&\vdots& \ddots
\end {array}$}}
\end{picture}
\end{center}
\caption{The $\SL{4}$-tiling associated to $\bleu{\ldots xxyyxxyyxxyy \ldots}$}\label{Fig_mot_slk}
\end{figure}
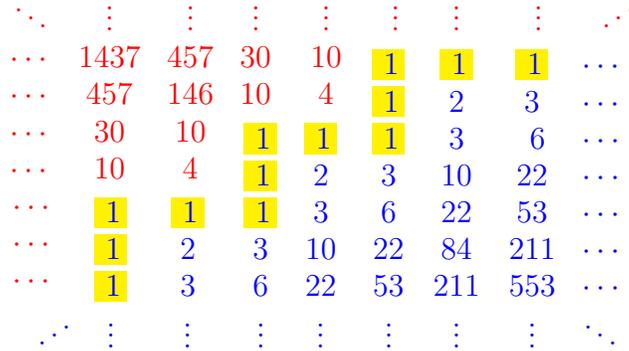

Let $p$ be a point lying below the path $\pi_w$, and suppose that $w_p$ factors as $x^i\,u\,y^j$ (with $i$ and $j$ maximal). Then we say that $u_p:=u$ is the \defn{short projection word} of $p$ on $w$. 
Illustrating with the tiling of Figure~\ref{Fig_mot_slk}, one may check that for the points $p$ corresponding to the entries with value equal to $6$ (lying below the path), one has $w_p=xxyyxxyy$ and $u_p=yyxx$ (for all instances of $6$); whereas for the points $p$ corresponding to the entries $30$ (lying above the path), one has $w_p= \bar {y}\bar{y}\bar{x}\bar{x}$ (likewise for all instances of $30$).


For a point $p$ lying below a path $\pi$, the \defn{distance} between $p$ and  $\pi$ is  the unique integer $k\in\N$ such that $(i-k+1,j-k+1)$ lies on $\pi$. Observe that there is but one point of $\pi$ lying on any given diagonal $x=y+c$. Our definition makes it so that points lying on the path are considered to be at distance $1$ of it (this will make our life easier later). 
We further consider the notion of \defn{$k$-fringe}, $\Phi_k(\pi)$ of a path $\pi$, i.e.: the points lying below the path that are within distance $k$ of it. Thus, we have
   \begin{equation}\label{def_fringe}
      \Phi_k(\pi):=\{(i+m,j+m)\ |\ (i,j)\in\pi,\quad{\rm and}\quad 0\leq m<k\}
   \end{equation}
Some of these notions are illustrated in Figure~\ref{fig_dessous}.

Given two points $p$ and $q$ on the path $\pi_w$ of a word $w$, we consider the sets of paths
\begin{equation}\label{defn_ense_chemins}
   \bleu{\Paths_{w;k}(p,q):=\{\theta\ {\rm a\ path}\ |\ \debut(\theta)=p,\  \fin(\theta)=q,\  {\rm and} \  \theta\subset \Phi_k(\pi_w)\}},
 \end{equation}
that start at $p$, end at $q$, while staying inside the $k$-fringe of $\pi_w$. We also denote by $\Paths_{w;k}$ the set of all such paths, disregarding start and end points.
The tiling $\A_{w;k}=\left(a_{ij}\right)_{i,j}$ is then defined, for points $p=(i,j)$ lying below the path $\pi_w$, by setting
\begin{equation}\label{def_tiling}
    \bleu{a_{ij}:=\# {\Paths_{w;k}(\chi_w(p),\gamma_w(p))}}.
    \end{equation}
For instance, for the word $w=\ldots yyxxyxyyyx\ldots $ and $k=2$, the resulting (partial) tiling is as follows:
$$\\A_{w;2}=\begin {array}{ccccccc} 
&&&&&&\r{1}\\ 
&&&\r{1}&\r{1}&\r{1}&\r{1}\\ 
&&\r{1}&\r{1}&2&3&4\\ 
&&\r{1}&2&5&8&11\\ 
\r{1}\ &\r{1}\ &\r{1}&3&8&13&18\end {array}$$
Observe that, for fixed $d$, any path $\theta$ in the $k$-fringe of $\pi$ goes through at most one of the the $k$ points of the set $\Phi_k(\pi)\cap \triangle_d$, where  $\triangle_d$ denotes the \defn {diagonal} 
\begin{equation}
     \triangle_d:=\{\ (i,j)\in \Z\times \Z \ |\ i-j=d\ \}.
\end{equation}
It is useful to have the following terminology: given a bi-infinite path and a tiling $\A$, let us call \defn {principal minors of order}  $m$ (relative to $\pi_w$), the minors of $\A$ that are of the form $M_{ij}^{(m)}$, with $(i,j)$ lying on the path $\pi_w$. In other word, a principal minor is an adjacent minor of $\A$ whose upper left corner lies on the path (in particular, it is contained in the $k-1$ fringe of the path when $m\leq k$). 

For $h\in\Z$, we also say that a minor $M_{ij}^{(m)}$ is \defn{located on the $h$-th diagonal} if we have $h=j-i$.
To tie all this to our study of $\SL{k}$-tilings, we now give entirely combinatorial arguments for the following statements.

\begin{proposition}\label{prop_word_tiling}
 \bleu{The partial tiling $\A_{w;k}$ is a $0$-free $\SL{k}$-tiling with principal minors of order $<k$ all equal to $1$.}
    \bleu{Moreover, for any $k$-subsets $I$ and $J$ such that $I\times J$ is contained in the shape of $\A$ {\rm (}all the points lying below $\pi_w${\rm )},
         we have}
\begin{equation}\label{mineurs_multiplicatifs_chemins}
    \bleu{\M_{IJ}=\M_{I,\{j,\ldots,j+k-1\}} \, \M_{\{i,\ldots,i+k-1\},J}}.
  \end{equation}    
\bleu{Finally, if the path $\pi_w$ is admissible, then $\A_{w;k}$ extends uniquely to a complete tame $\SL{k}$-tiling.}
    \end{proposition}

There is some redundance here, since the 0-free property implies that $\A$ is tame (Proposition~\ref{rankk}), and thus  (\ref{mineurs_multiplicatifs_chemins}) holds by Proposition~\ref{rank1}. Notwithstanding, we want to make evident that nice combinatorial methods may be used to understand all this. 

\begin{proof}[\bf Proof.]
It follows from a theorem of Gessel-Viennot (see \cite{viennot}) that 
 we may interpret combinatorially  any minor of  $\A_{w;k}$ as follows. 
 Recall that a family of paths is said to be non-crossing if no pair of paths in the family has a common point.
Given equal cardinality subsets $I$ of rows and $J$ of columns, such that $I\times J$ lies below $\pi_w$, we denote by $\NC_{IJ}$ the set of all non-crossing families of paths in $\Paths_{w;k}$ linking the horizontal projection of $I$ on $\pi_w$ to the vertical projection of $J$ on $\pi_w$. More precisely, let
    $$I=\{i_1<i_2<\ldots <i_m\}, \qquad{\rm and}\qquad J=\{j_1<j_2<\ldots <j_m\},$$
and denote by $p_1,\ldots, p_m$ and $q_1,\ldots, q_m$ the respective horizontal and vertical projections on $\pi_w$.
Then the elements  of $\NC_{IJ}$ are ``configurations'' $\{\pi_1,\pi_2,\ldots,\pi_m\}$  of paths $\pi_s $ in $\Paths_{w;k}$, with  
      \begin{equation}
             \bleu{ \debut(\pi_s)=p_s}, \quad {\rm and}\quad \bleu{\fin(\pi_s)=q_s},
       \end{equation}
no two of which cross.  In our context, the aforementioned theorem of \cite{viennot} states that we have
\begin{equation}\label{mineurs_chemins}
    \bleu{\M_{IJ} = \# \NC_{IJ}}.
  \end{equation}
Recall that this is shown by constructing a sign changing involution on the set of crossing configurations, thus showing that they can be eliminated from a global signed counting that clearly corresponds to the evaluation of the determinant considered.

Observe that non-crossing path configurations $\{\pi_1,\pi_2,\ldots,\pi_m\}$ intersect any given diagonal $\Delta_d$ in at most $m$ points. In fact, this intersection number is exactly equal to $m$ for the diagonals that pass through points of $\pi_w$ lying between $p_1$ and $q_1$. This forces
all the sets $ \NC_{IJ}$ to be empty whenever $\#I=\#J>k$. Hence the corresponding minors all vanish, so that the tiling is of rank $k$. 

To continue with our combinatorial argument, let us write $ \NC_{ij}^{(m)}$
when 
    $$I=\{i,\ldots,i+m-1\}, \quad {\rm and}\quad J=\{j,\ldots,j+m-1\},$$  
 so that we have
\begin{equation}\label{mineurs_chemins_ij}
    \bleu{\M_{ij}^{(m)} = \# \NC_{ij}^{(m)}}.
  \end{equation}
In the $k$-fringe of $\pi_w$, there is room for one and exactly one configuration of $k$ non-crossing paths having adjacent starting points and adjacent end points, so that
we must necessarily have $\M_{ij}^{(k)} =1$ for all point $(i,j)$ in the tiling, hence the $\SL{k}$-condition is verified.

The tiling is $0$-free, and in fact we have $\M_{ij}^{(m)}\geq 1$ for all $m$ between $1$ and $k$, since there exists corresponding non-crossing path configurations for all these $m$.
Finally, the multiplicative property (\ref{mineurs_multiplicatifs_chemins}) of $k\times k$-minors can easily be explained as follows, in terms of configurations of $k$ non-crossing paths.
With $I$ and $J$ satisfying the hypothesis of the proposition, the required identity  follows from a simple bijection
  \begin{equation}\label{chemins_multiplicatifs}
    \bleu{\NC_{IJ} \rightarrow \NC_{I,\{j,\ldots,j+k-1\}} \times \NC_{\{i,\ldots,i+k-1\},J}}, 
   \end{equation}
obtained by breaking up the paths considered into three portions as follows. For a path $\pi_s$ (starting at $p_s$ and ending at $q_s$) in $\Phi_k(\pi_w)$, denote respectively by $p$ and $q$ the points of $\pi_s$ that lie on the diagonals respectively containing $p_1$ and $q_1$. These exist since $p_s\leq p_1\leq q_1\leq q_s$. We decompose $\pi_s$ as the concatenation 
     $$\pi_s= \pi_s^{(1)} \pi_s^{(2)} \pi_s^{(3)},$$
 with 
 \begin{enumerate}
     \item[$\bullet$] $\pi_s^{(1)}$ being the portion of $\pi_s$ going from $p_s$ to $p$,  
     \item[$\bullet$] $\pi_s^{(2)}$ being the portion of $\pi_s$ going from $p$ to $q$,  and 
     \item[$\bullet$] $\pi_s^{(3)}$ being the portion of $\pi_s$ going from $q$ to $q_s$.
 \end{enumerate} 
In particular, all the paths $\pi_s^{(2)}$ start on the same diagonal (the one that contains $p_1$) and end on the same diagonal (the one that contains $q_1$). Since these $k$ paths are non-crossing and all lie in the $k$-fringe, there is but one possibility for the resulting configuration $\{\pi_1^{(2)},\ldots, \pi_k^{(2)}\}$. We easily identify the configurations $\{\pi_1^{(1)},\ldots, \pi_k^{(1)}\}$ with elements of $\NC_{I,\{j,\ldots,j+k-1\}}$ (by application of the same decomposition as above to these last elements, observing that in this cases third components are trivial). Likewise we identify the configurations $\{\pi_1^{(3)},\ldots, \pi_k^{(3)}\}$ with elements of $\NC_{\{i,\ldots,i+k-1\},J}$. This establishes the bijection.

The $\SL{k}$-tiling $\A_{w;k}$ may be uniquely completed into a tame $\SL{k}$-array by Lemma~\ref{lem_unique} below. 
 \end{proof}

For example, with $k=4$ and the word $w=\cdots xxyyxxyyxxyy\cdots$, we first get the partial array (lying below the path) of Figure~\ref{Fig_mot_slk} by path enumeration, and then complete it to get a $\SL{4}$-tiling of $\Z\times \Z$. 

\begin{lemma}\label{lem_unique}
\bleu {Consider a partial $\SL{k}$-tiling which is defined on every point below a given path, and such that all of its adjacent $(k+1)\times (k+1)$-subminor  (lying entirely in its shape)  vanishes. If the path is admissible, then the partial $\SL{k}$-tiling extends uniquely to a complete tame $\SL{k}$-tiling.}
\end{lemma}

\begin{proof}[\bf Proof.]
Let $i$, $j$ and $k$ be integers such that $\{i,i+1,i+2,\ldots\} \times \{j-k,\ldots,j-1,j\}$ is contained in the shape of $\A$. denote by $C_{j-k},\ldots,C_{j-1},C_j$ the columns of the corresponding submatrix of $\A$. Then, it follows from the vanishing $(k+1)\times (k+1)$-subminor condition, that we have a relation of the form
 $$\bleu{ (-1)^kC_{j-k}-(-1)^{k}a_1C_{j-k+1}+\cdots -a_{k-1}C_{j-1}+C_j=0}.$$
 for some $a_1,\ldots, a_{k-1}$ in $\K$. Note that the coefficients $a_h$ are independent of the  $i$ chosen.

Considering the analogous argument for rows, and assuming that the origin of the plane is in the shape of $\A$, we obtain a linearization data (see Section~\ref{sec_general}). Using this linearization data, we may apply Proposition~\ref{coefflinear} to get a complete tame tiling of the plane. Call it $\A'$. It follows from the construction that $\A$ and $\A'$ coincide on the shape of $\A$, which proves the lemma.
\end{proof}

\subsection*{Weighted word tilings} We now extend the previous construction to the situation where paths are given Laurent monomial weights. 

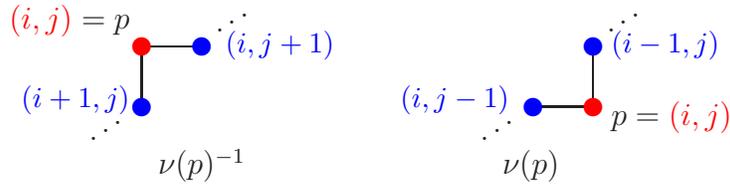
\begin{figure}[ht]\setlength{\unitlength}{8mm}
 $$ \begin{matrix}
    \begin{picture}(4,2)(-1,-.5)
  \anticoude(0,0)11
   \D(0,0){1}\X(0,1){0}{1}\D(1,1){1} 
   \put(1.2,1.2){$\revddots$}
   \put(-0.9,-0.6){$\revddots$}
   \put(1.4,.9){\small$\bleu{(i,j+1)}$}
    \put(-2,0){\small$\bleu{(i+1,j)}$}
    \put(-2.2,1.3){$\rouge{(i,j)}=p$} 
   \end{picture}&\qquad\qquad
   \begin{picture}(4,2)(-1,-.5)
  \coude(0,0)11
   \D(0,0){1}\X(1,0){0}{1}\D(1,1){1}
   \put(1.2,1.2){$\revddots$}
   \put(-0.9,-0.6){$\revddots$}
   \put(1.3,.9){\small$\bleu{(i-1,j)}$}
   \put(-2.2,0){\small $\bleu{(i,j-1)}$}
   \put(1.3,-0.3){$p=\rouge{(i,j)}$}
  \end{picture}\\
   \nu(p)^{-1}& \nu(p)
\end{matrix}$$
\caption{Right and left turns at $p$, and corresponding weight.}\label{Fig_turn}
 \end{figure}

At a point $p=(i,j)$, along a path $\theta$, we say that we have a \defn{right-turn} (resp. \defn{left-turn}) if both  $(i+1,j)$ and $(i,j+1)$) (resp. $(i,j-1)$ and $(i-1,j)$) belong to the path $\theta$. This is illustrated in Figure~\ref{Fig_turn}.
 For a given (bi-infinite) word $w$,  we start by giving a weight $\nu(p)$ to each point $p=(i,j)$ in the $k$-fringe of $\pi_w$, setting
     $$\nu(p):=  \frac{t_{j-i,r}}{t_{j-i,r-1}},$$
where $r$ is the distance between $p$ and $\pi_w$, We assume here that the $t_{m,r}$ are independent commuting variables, setting $t_{m,r}=1$ whenever $r\leq 0$ or $r\geq k$.
With all this at hand, define the \defn{weight}  $\omega(\theta)$ (with respect to the word $w$) of  a nonempty path $\theta$ to be the product 
    \begin{equation}\label{defn_poids_chemin}
        \bleu{\omega(\theta):=\alpha\, \beta \prod_{p\ {\rm left-turn\ of\ }\theta} \nu(p)\ \cdot \prod_{p\ {\rm right-turn\ of\ }\theta} \nu(p)^{-1}},
     \end{equation}
where we set $\alpha:= \nu(p_s)$, if $\theta$ starts at $p_s$  by a vertical step. Otherwise we set $\alpha=1$. Likewise, we set $\beta:=\nu(p_e)$, if  $\theta$ ends at $p_e$ after an horizontal step. Otherwise, we set $\beta=1$. Finally, when $\theta $ is the empty path, both starting and ending at $p$, we  simply set $\omega(\theta):=\nu(p)$.

We then consider the partial tiling $\mathcal{B}_{w;k}:=(b_{ij})_{i,j}$, for point $(i,j)$ lying below the path, obtained by setting
 \begin{equation}\label{pavage_pondere}
      \bleu{b_{ij}:=\sum_{\theta} \omega(\theta)},
 \end{equation}
 for $\theta$ varying in the set $\Paths_{w;k}(\chi_w(i,j),\gamma_w(i,j))$, of paths starting at $\chi_w(i,j)$ and ending at $\gamma_w(i,j)$.    
 

 \begin{proposition}  \label{pavage_variables_t}
 \bleu{There is a unique tame $\SL{k}$-tiling of $\Z\times \Z$ extending $\mathcal{B}_{w;k}$, with entries Laurent polynomials in the variables $t_{hr}$. More precisely, the values are in the subsemiring generated by these variables and their inverses. Moreover, each principal minor of order $r$, $r=1,\ldots,k-1$, located on the $h$-th diagonal, is equal to $t_{hr}$.}
 \end{proposition}
 
\begin{proof}[\bf Proof.]
Again, we simply apply the Gessel-Viennot technique, verifying that the involution (as in their original proof), required to show that crossing path configurations may be eliminated, is weight preserving. There are several cases, left to the reader.
It follows, as in the proof of Proposition~\ref{prop_word_tiling}, that $b_{ij}$ is a tame $\SL{k}$-tiling.  It is clear that (\ref{pavage_pondere})   is in the described semiring.
Moreover, by the non-crossing path description, each principal minor of order $r<k$ is equal to $t_{hr}$, if the minor is located on the $h$-th diagonal.
\end{proof}

This proposition may be used to construct $\SL{k}$-tilings having arbitrary values (variables) as entries in the $(k-1)$-fringe of the path $\pi_w$. It turns out that the entries of the resulting tiling actually lie in the polynomial ring generated by these variables as well as the inverses of all principal minors (relative to $\pi_w$) of order at most $k-1$. This is a Laurent-like phenomenon (see \cite{fomin}) of a new kind. Moreover, we may in fact replace ``ring''  by ``semiring'', so that we actually get a positivity result, just as is the case in the theory of cluster algebras.

\begin{corollary}
\bleu {An admissible path $\pi_w$ being given, associate to each point in its $(k-1)$-fringe, a distinct commutative variable. Then this assigment extends uniquely into a complete tame $\SL{k}$-tiling of the plane whose values are in the semiring generated by the principal minors of order $<k$ and their inverses}.
\end{corollary}

\begin{proof}[\bf Proof.]
Consider the tiling of Proposition~\ref{pavage_variables_t}. Let $a_{ij}$ denote its value at the point $(i,j)$. Then, by the same proposition, each $a_{ij}$ is in the semiring generated by the variables $t_{hr}$ and their inverses, where $h\in \Z$ and $r=1,\ldots,k-1$. 

Let $s_{hr}=t_{hr}/t_{h,r-1}$ for $h\in \Z$ and $r=1,\ldots,k-1$. Recall that the $t_{hr}$ are distinct commuting variables, and that $t_{h,0}=1$. The field of fractions $\K$ in the variables $t_{hr}$ is also generated by the $s_{hr}$, and the mapping $t_{hr}\mapsto s_{hr}$ is an automorphism of this field. 
For $h\in \Z$ and $r$ going from $1$ to $k-1$, denote by $\alpha_{hr}$ the entry $a_{ij}$, if the point $(i,j)$ lies at distance $r$ below the path $\pi_w$, on the $h$-th diagonal. 
By the path description, we see that, $\alpha_{hr}$ is the sum of $s_{hr}$ and of a fraction in the $t_{h'r'}$ with $(h',r')<(h,r)$ for the natural order on $\Z^2$. The latter fraction, when expressed in the $s_{h'r'}$, involves only variables $s_{h'r'}$ with the same condition. Hence, the function $t_{hr}\mapsto \alpha_{hr}$ defines an automorphism of $\K$.

Now, let $x_{hr}$ be a family of distinct commuting variables, for $h\in \Z$ and $r=1,\ldots,k-1$. Let $\L$ be its field of fractions. The fields $\K$ and $\L$ are isomorphic (e.g. by the mapping $t_{ij}\mapsto  x_{ij}$). By what we have just seen, the mapping $\alpha_{hr}\mapsto x_{hr}$, $h\in \Z$ and $r=1,\ldots,k-1$ defines an isomorphism from $\K$ onto $\L$. If we map each $a_{ij}$ under this isomorphism, we obtain a tame $\SL{k}$-tiling $X=(b_{ij})$ such that $b_{hr}=x_{hr}$ for $h\in \Z$ and $r=1,\ldots,k-1$. This also implies that we may find elements $\tau_{hr}$ in the field $\L$ such that $b_{ij}$ is in the semiring generated by the $\tau_{hr}$ and their inverses, $h\in \Z$ and $r=1,\ldots,k-1$.
Furthermore, by Proposition~\ref{pavage_variables_t}, the principal $r\times r$-minor of $X$ ($r<k$), located on the $h$-th diagonal, is equal to $\tau_{hr}$.

Unicity follows from the following lemma, which of independent interest.
\end{proof}

\begin{lemma}
\bleu{An admissible path $\pi_w$ being given, associate to each point in its $(k-1)$-fringe, an element of some field. Suppose that the $(k-1)$-principal minors relative to $\pi_w$ are nonzero. Then this partial tiling extends uniquely to a tame $\SL{k}$-tiling of the plane.}
\end{lemma}

\begin{proof}[\bf Proof.]
Indeed, under the nonzero $(k-1)$-principal minor hypothesis, the $\SL{k}$ property imposes that we have a unique extension of the partial tiling to its $k$-fringe. This furnishes enough $k\times (k+1)$ and $(k+1)\times k$ submatrices so that we may compute the linearization data for any tame $\SL{k}$-tiling that would extend the $k$-fringe shaped partial tiling (see the remark following the proof of Proposition~\ref{coefflinear}). This proves unicity, in view of the same proposition.
Existence, which will not be used here, is left to the reader.
\end{proof}

\begin{figure}[ht]\setlength{\unitlength}{6mm}
\bleu{\begin{picture}(8,8)(0,0)
  \put(2,8){ \put(0,0){$1$}\put(1,0){$1$}\put(2,0){$1$}\put(3,0){$1$}\put(4,0){$1$}\put(5,0){$1$}}
  \put(1,6){ \put(0,0){$1$}\put(1,0){$1$}\put(2,0){$2$}\put(3,0){$3$}\put(4,0){$4$}\put(5,0){$\cdots$}}
  \put(1,4){ \put(0,0){$1$}\put(1,0){$2$}\put(2,0){$5$}\put(3,0){$9$}\put(4,0){$\cdots$}}
  \put(0,2){ \put(0,0){$1$}\put(1,0){$1$}\put(2,0){$3$}\put(3,0){$9$}\put(4,0){$\cdots$}}
    \put(0,0){$1$}\put(1,0){$2$}\put(2,0){$7$}\put(3,0){$\cdots$}
 \end{picture}}
 \rouge{\begin{picture}(8,8)(7.8,-1)\tiny
   \put(2,6){\put(0,0){$1$} \put(1,0){$1$} \put(2,0){$1$}\put(3,0){$1$}\put(4,0){$1$}}
   \put(1,4){\put(0,0){$1$} \put(1,0){$1$} \put(2,0){$3$} \put(3,0){$6$}\put(4,0){$\cdots$}}
   \put(1,2){\put(0,0){$1$} \put(1,0){$3$} \put(2,0){$14$} \put(3,0){$\cdots$}}
   \put(0,0){$1$} \put(1,0){$1$} \put(2,0){$6$}\put(3,0){$\cdots$}
 \end{picture}}
\caption{Joint calculation of a $\SL{3}$-tiling and its dual.}\label{pavage_et_dual}
\end{figure}
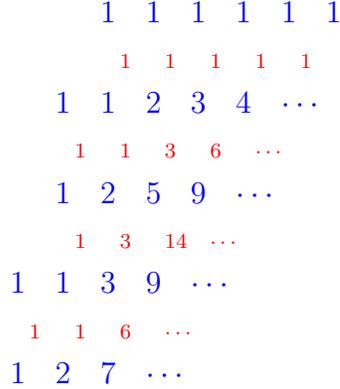

In \cite{assem}, one may find many $\SL{2}$-tilings associated to paths, both over the integers, and with arbitrary variables on the path. 
The case $\SL{3}$ has an extra interesting feature. Indeed, a consequence of (\ref{rec_trois}) is that we can very elegantly characterize any $\SL{3}$-tiling in tandem with its dual tiling. Indeed, under the assumption that $\A$ is $\SL{3}$ and writing $\A^*=\left(a^{*}_{ij}\right)_{i,j}$, the tiling identity considered is equivalent to the family of  equalities
  \begin{equation}\label{outils_SL3}
      \begin{array}{rcl}
        \bleu{a_{ij}}&=&\bleu{\displaystyle \frac{1}{a_{i-1,j-1}}(\rouge{a^{*}_{i-1,j-1}} + a_{i-1,j}a_{i,j-1})},\\[12pt]
        \rouge{a^{*}_{ij}}&=&\bleu{\displaystyle\frac{1}{\rouge{a^{*}_{i-1,j-1}}}(a_{i-1,j-1} + \rouge{a^{*}_{i-1,j}a^{*}_{i,j-1}})}.
    \end{array}
  \end{equation}
 This makes it evident (in another fashion) that the tiling constructed from a path is positive (and non-zero) for points lying below the path, since entries of $\A$ and $\A^*$ may be calculated recursively in parallel  using the positive expression on the right-hand-side of (\ref{outils_SL3}).  In the case of integer tilings, this is illustrated in Figure~\ref{pavage_et_dual}. Large entries correspond to the $\bleu{a_{ij}}$'s, and smaller ones correspond to the $\rouge{a^*_{ij}}$'s. The entry $\rouge{a^*_{ij}}$ sits immediately to the south-east of $\bleu{a_{ij}}$. Clearly the recursion process may be continued where it is left off. It corresponds to the statement that each number is obtained as the determinant of the $4$ numbers that immediately surround it.

\section{Matrix description}
Consider the morphism $\mu$, from the free group $F_{x,y}$ (on the letters $x$ and $y$) to the group $\SL{k}$, which is obtained by setting
\begin{equation}
 \bleu{\mu(x):=\Id + N}, \qquad{\rm and}\qquad
       \bleu{\mu(y):=\Id+N^{\rm tr}},
   \end{equation}
 where we denote by $N$ the matrix nilpotent $k\times k$ matrix
    $$N:=\begin{pmatrix} 
    0& 1& 0   &&\ldots &0\\
      & 0 & 1 &0 & &         \\ 
    & & \ddots   &\ddots&\ddots &\vdots \\ 
        &  & &0&1 &0\\
    0&   &\ldots & &0 &1
    \end{pmatrix}$$
(Recall that $N$ is nilpotent of order $k$, so that $N^k=0$, and $N^i\not=0$ when $i<k$.)
    We denote $\overline{x}$ the inverse of $x$ in the free group $F_{x,y}$, and likewise for $y$. 
    We then define the function ${\mathcal T}_w:\Z\times \Z\rightarrow \N$ as
\begin{equation}\label{pavage_mot_matrice}
\bleu{{\mathcal T}_w(p):= 
        {\bf e}_k\,  \mu(w_p)\, {\bf e}_k^{\rm tr}}         
        \end{equation}
where ${\bf e}_k$ denotes the unit $k$-vector $(0,\ldots,0,1)$. Recall that the projection word $w_p$ has been defined in Section~\ref{words_tilings}.

\begin{proposition}\label{theo_neuf}
\bleu{For any admissible bi-infinite word $w$, the function ${\mathcal T}_w$ is a tame $\SL{k}$-tiling, whose principal minors of order $<k$ are all equal to $1$. It coincides with the tiling of Proposition~\ref{prop_word_tiling}.}
\end{proposition}

For the proof of Proposition~\ref{theo_neuf}, see Section~\ref{preuves}. Observe that this result easily implies the following:

\begin{corollary}
\bleu{With the same hypothesis as in Proposition~\ref{theo_neuf}, we have}
\begin{itemize}
   \item[(i)] \bleu{If $p$ lies below $\pi_w$ then} 
    \begin{equation}
        \bleu{ {\mathcal T}_w(p)= {\bf e}_k\,  \mu(u_p)\, {\bf e}_k^{\rm tr}}.
    \end{equation}
\bleu{In other words we can replace the projection word $w_p$ by the short projection word $u_p$ in our calculations}.
 \item[(ii)] \bleu{If $p$ lies above $\pi_w$,with $w_p=\overline{x_1}\cdots \overline{x_n}$, $x_1,\ldots,x_n \in \{x,y\}$, then}
   \begin{equation}
       \bleu{ {\mathcal T}_w(p)={\bf e}_k\, \mu'({x_1}\cdots x_n)\, {\bf e}_k^{\rm tr}},
    \end{equation}
  \bleu{  where $\mu'$ is the morphism such that}
   $$\bleu{\mu'({x}):=(Id-N)^{-1}}, \qquad{\rm and}\qquad
       \bleu{ \mu'({y}):=\mu'({x})^{\rm tr}}.$$ 
\end{itemize}
\bleu{In particular, we conclude that ${\mathcal T}_w(p)$ is positive for all $p$.}
\end{corollary}

\begin{proof}[\bf Proof.] Assume that $w=w_p=x^i\,u\,y^j$ with $u=u_p$. To show (i),  we first observe that
\begin{eqnarray*}
     \mu(x^i)&=&(\Id+N)^i\\
                  &=& \Id+\sum_{j=0}^i \binom{i}{j} N^i
  \end{eqnarray*}
  is upper unitriangular, hence ${\bf e}_k \mu(x^i) = {\bf e}_k$.
 Likewise,  $\mu(y^i)$ is lower unitriangular, so that $\mu(y^j) {\bf e}_k^{\rm tr}= {\bf e}_k^{\rm tr}$. Thus we directly calculate that 
\begin{eqnarray*}
        {\mathcal T}_w(p) &=& {\bf e}_k\,  \mu(x^i u y^j)\, {\bf e}_k^{\rm tr}\\
                        &=&  {\bf e}_k\,  \mu(x^i)\mu(u)\mu( y^j)\, {\bf e}_k^{\rm tr}\\
                        &=&  {\bf e}_k\,  \mu(u)\, {\bf e}_k^{\rm tr}
 \end{eqnarray*}
as announced.

For (ii), we make use of the matrix isomorphism 
$$\alpha(A):=D_k\, A\, D_k^{-1},$$
with $D_k$ standing for the diagonal matrix with entry equal to $(-1)^{i+1}$ on the diagonal. Clearly, $\alpha(a_{ij})=((-1)^{i+j}a_{ij})$. Thus $\alpha(\mu'(x))=\mu(\overline{x})$ and $\alpha(\mu'(y))=\mu(\overline{y})$ as is easily verified. Hence, for any $x_1,\ldots,x_n \in \{x,y\}$, we have $\alpha(\mu(x_1\cdots x_n))=\mu(\overline{x_1}\cdots \overline{x_n})$. We conclude since $\alpha(A_{kk})=\A_{kk}$. For the final assertion, note that $\mu'$ has nonnegative coefficients.

\end{proof}

\section{Proofs}\label{preuves}
To prove some of our previous assertions we first need a few 
linear algebra lemmas.

\begin{lemma}\label{det_rank}
\bleu{If a matrix has all its adjacent $(k+1)\times (k+1)$ minors vanishing, whereas no adjacent $k\times k$ vanishes, then it is of rank $k$}.
\end{lemma}

\begin{proof}[\bf Proof.]
It is enough to show that for any choice of $k+1$ successive columns $C_0,\ldots, C_k$ of this matrix, $C_0$ (resp. $C_k$) is a linear combination of $C_1, \ldots, C_k$ (resp. of $C_0, \ldots, C_{k-1}$). By symmetry, the property with $C_0$ will suffice.  Let $v_i$ denote the rows of the matrix $(C_0,\ldots,C_k)$. Note that $v_i$ is of length $k+1$.

To show our assertion, let us construct a non-vanishing linear form $\varphi$ that annihilates all $v_i$. The existence of such a linear form implies the existence of scalars $a_0,\ldots, a_{k}$ such that $\sum_{j=0,\ldots,k} a_jC_j=0$. Moreover $a_0$ has to be nonzero, since otherwise it would contradict the assumption on the non-vanishing $k$-minors. 

Such a linear form exists for $k+1$ successive rows of $M$, since $\det(M)=0$ by assumption. Consider $k+2$ successive rows of $M$, and two non-vanishing linear forms $\varphi$ and $\psi$ such that the first $k+1$ rows are in $\Ker(\varphi)$ and the $k+1$ last are in $\Ker(\psi)$. Then we argue as follows to show that $\varphi$ and $\psi$ must necessarily be proportional. If we restrict the two linear forms to the $k$ intermediate rows, $v_1,\ldots, v_k$ say, we see that $\varphi$ and $\psi$, considered as column vectors of length $k+1$, are both annihilated by the $(k\times (k+1)$-matrix
\begin{equation}
\begin{pmatrix}
v_1\\v_2\\\vdots\\v_k
\end{pmatrix},
\end{equation}
whose rows are the vectors $v_i$.
By assumption, this matrix is of rank $k$, hence it has a kernel of dimension $1$. It follows that its columns vectors are proportional, and thus so are $\varphi$ and $\psi$.
\end{proof}

\begin{lemma}\label{rang1}
\bleu{Let $A$ be a square matrix of order $k+1$ such the matrix of its $k\times k$-minors $(\det(A_{IJ})_{I,J}$, with $I$ and $J$ running through all $k$-subsets of  $\{1,\ldots,k+1\}$, is of rank $1$. Then $\det(A)=0$}.
\end{lemma}

\begin{proof}[\bf Proof.]
If the central $(k-1)\times (k-1)$-minor $\det A_{\{2,\ldots,k\},\{2,\ldots,k\}}$ of $A$ is nonzero, then  (\ref{carrol}), with $r=k$,  implies that $\det(A)=0$. If some $(k-1)\times (k-1)$-minor of $M$ is nonzero, we may bring it into central position by row and column permutations; 
these operations amount to row and column permutations of the matrix of $k\times k$-minors of $A$; hence, by the previous argument, $\det(A)=0$. Finally, if all the $(k-1)\times (k-1)$-minors of $A$ vanish, then so does $\det(A)$.
\end{proof}

\begin{lemma}\label{lemme_1}
\bleu{Let $V$ be a vector space, and consider a finite ordered set of indices $K$ for which we have selected vectors
$v_k,v'_k$ in $V$, as well as   $u_k, u'_k$ in the dual space $V^*$. Assuming that for all $k$ in $K$ we have the relations\footnote{Here, as in the sequel of this section, the stars {\rm (}$\etoile${\rm )} stand for some coefficients that we do not actually need to specify.}
   \begin{eqnarray*}
           v'_k&=& v_k +\sum_{\ell<k} \etoile\, v_\ell,\\
           u'_k&=& u_k +\sum_{\ell<k} \etoile\, u_\ell.
    \end{eqnarray*}
Then we have the equality}
  \begin{equation}
    \bleu{ \det\left(u_k(v_\ell)\right)_{k,\ell\in K} =\det \left(u'_k(v'_\ell)\right)_{k,\ell\in K}} .
  \end{equation}
\end{lemma}

\begin{proof}[\bf Proof.]
We simply pass from one matrix to the other by multiplication on the left and on the right by uni-triangular matrices.
\end{proof}

Consider now 
intervals of cardinality $k-1$ of the set $\{2,\ldots, 2k-1\}$, of the form
    $$I_q:=\{q+1,q+2,\ldots, q+k-1\},\quad q=1,\ldots, k.$$
 For convenience sake, we write $K_q:=[k]\setminus \{q\}$ (with $[k]$ standing as usual for $\{1,\ldots,k\}$). Let $e_1,\ldots, e_{2\,k-1}$ be elements of some vector space $V$. For $J=\{j_1,\ldots,j_s\}$ such that
         $$1\leq j_1\leq \ldots\leq j_s\leq 2k-1,$$
we denote by $e_J$ the wedge product
        $$e_J:=e_{j_1}\wedge e_{j_2}\wedge\ldots\wedge e_{j_s}.$$ 
Then the following holds. 

\begin{lemma}\label{lemme_2}
\bleu{   For vectors  $e_1,\ldots, e_{2\,k-1}$ in $V$ which are such that}
         $$\bleu{e_j=(-1)^{k-1}e_{j-k}+\etoile\, e_{j-k+1} +\ldots + \etoile\, e_{j-1},\qquad {\rm for}\quad j=k+1,\ldots, 2k-1},$$
   \bleu{the following identity holds for all $q=1,\ldots, k$}:
        \begin{equation}\label{formule_lemme_2}
            \bleu{ e_{I_q}= (-1)^{q+1} e_{K_q}+ \etoile\, e_{K_{q-1}}+\ldots+ \etoile\, e_{K_1}}.
        \end{equation}
\end{lemma}

\begin{proof}[\bf Proof.]
Writing $E:=e_{q+1}\wedge e_{q+2}\wedge\ldots\wedge e_k$ and $\varepsilon:=(-1)^{k-1}$, we calculate that
\begin{eqnarray*}
   e_{I_q}&=& E\wedge e_{k+1}\wedge \ldots\wedge e_{k+q-1}\\
              &=& E\wedge e_{k+1}\wedge \ldots\wedge e_{k+q-2}\wedge(\varepsilon e_{q-1}+\etoile\,e_q+\ldots + \etoile\, e_{k+q-2} ) \\
              &=& E\wedge e_{k+1}\wedge \ldots\wedge e_{k+q-2}\wedge(\varepsilon e_{q-1}+\etoile\,e_q ) \\
              &&\qquad\hbox{\rm (since $e_{q+1},\ldots e_{k+q-2}$ appear as factors in the product $E\wedge e_{k+1}\wedge \ldots \wedge e_{k+q-2}$.)}\\
               &=& E\wedge e_{k+1}\wedge \ldots\wedge e_{k+q-3}\wedge(\varepsilon e_{q-2}+\etoile\, e_{q-1} +\ldots + \etoile\, e_{k+q-3})\\ 
               &&\hskip1.7in            \wedge(\varepsilon e_{q-1}+\etoile\,e_q ) \\
               &=& E\wedge e_{k+1}\wedge \ldots\wedge e_{k+q-3}\wedge(\varepsilon e_{q-2}+\etoile\, e_{q-1} +\etoile\, e_{q})\\ 
               &&\hskip1.7in            \wedge(\varepsilon e_{q-1}+\etoile\,e_q ) \\
               &=&\ldots\\
               &=&E\wedge (\varepsilon e_1+\etoile\, e_2+\ldots+\etoile\,e_q) \wedge\ldots   \\
               &&\qquad          \wedge(\varepsilon e_{q-2}+\etoile\, e_{q-1} +\etoile\, e_{q})\\ 
               &&\qquad           \wedge(\varepsilon e_{q-1}+\etoile\,e_q ) \\
               &=&(-1)^{(q-1)(k-q)} (\varepsilon e_1+\etoile\, e_2+\ldots+\etoile\,e_q)  \wedge\ldots \wedge(\varepsilon e_{q-2}+\etoile\, e_{q-1} +\etoile\, e_{q})\\ 
               &&\qquad           \wedge(\varepsilon e_{q-1}+\etoile\,e_q )\wedge E. 
\end{eqnarray*}
The product that precedes $E$ is evidently in the $(q-1)^{\rm th}$-exterior power of the span of $e_1,\ldots,e_q$. It is thus a linear combination of the $e_{[q]\setminus \{i\}}$, for $i=1,\ldots,q$. It follows (as we are multiplying these $e_{[q]\setminus \{i\}}$ by $E$ on the right) that we have expressed $e_{I_q}$ as a linear combination of the $E_{K_i}$, for $i=1,\ldots, q$. Moreover, $e_{K_q}$ appears only once in the resulting expression. Its coefficient is thus $(-1)^{(q-1)(k-q)}(-1)^{(q-1)(k-1)}$. We conclude that (\ref{formule_lemme_2}) holds, since   $(q-1)(k-q)+(q-1)(k-1)\equiv (q-1)(-q-1)\equiv (q+1)^2\equiv(q+1)$ modulo $2$.
\end{proof}
In the next result, the first row and first column of matrices are indexed by $1$.

\begin{proposition}\label{Condense}
\bleu{Consider a  $(2k-1)\times(2k-1)$ matrix $A$ of rank $k$ having all of its adjacent $k\times k$ minors equal to $1$, and write $B=A_{11}^{(k)}$, $C=\partial_{k-1}A$, and $D=C_{22}^{(k)}$. Then, for all $h\leq k$, we have
    $$\det D_{11}^{(h)}=\det B_{h+1,h+1}^{(k-h)}.$$}
\end{proposition}

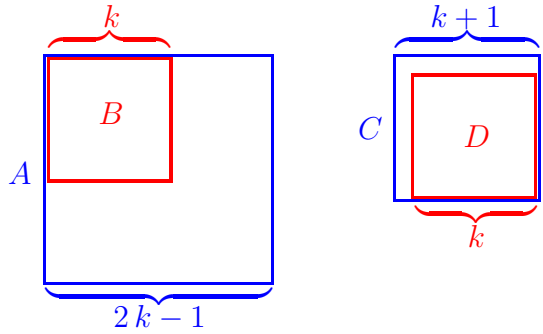
\begin{figure}[ht]
\setlength{\unitlength}{2mm}
    \begin{center} \label{dessin}
        \begin{picture}(15,19)(-3,-1.5)
                   \put(3.5,10.5){$\rouge{B}$}
                   \put(3.8,17){$\rouge{k}$}
                  \put(0.2,15.5){\rouge{$\overbrace{\hskip1.6cm}$}}
                   \put(-2.5,6.5){$\bleu{A}$}
                   \put(0.25,6.8){\rouge{\carre(0,0)(8)}}
                   \bleu{\carre(0,0)(15)}
                   \put(4.5,-3){$\bleu{2\,k-1}$}
                   \put(0,-0.3){\bleu{$\underbrace{\hskip3cm}$}}
                  \end{picture}\qquad \qquad
        \begin{picture}(15,16)(-3,-1)
                   \bleu{\carre(0,6)(9.5)}
                   \put(1.2,6.2){\rouge{\carre(0,0)(8)}}
                   \put(4.5,9.5){$\rouge{D}$}
                   \put(4.8,2.8){$\rouge{k}$}
                   \put(1.3,5.7){\rouge{$\underbrace{\hskip1.6cm}$}}
                   \put(2.3,17.5){$\bleu{k+1}$}
                  \put(0,16){\bleu{$\overbrace{\hskip 1.9cm}$}}
                   \put(-2.5,10){$\bleu{C}$}
                  \end{picture}
      \caption{The square matrices $A$, $B$, $C$, and $D$ of Proposition~\ref{Condense}.} \label{fig_ABCD}
    \end{center}
 \end{figure}

Observe that the square matrices $A,B,C,D$ are respectively of order $2k-1$, $k$, $k+1$ and $k$ (as illustrated in Figure~\ref{fig_ABCD}). Recall also  that $\partial_{k-1}A$ is the matrix of adjacent $k-1$-minors of $A$.

\begin{proof}[\bf Proof.]
Let $e_1,\ldots,e_{2k-1}$ be the column vectors of $A$, and consider the vector space $V$ that they span. Dually, let $\varphi_1,\ldots,\varphi_{2k-1}$ be the restriction to $V$ of the 
$2k-1$ projections of column vectors on the underlying field of scalars.

We clearly have $\varphi_i(e_j)=a_{ij}$. Using the usual duality\footnote{$\langle \psi_1\wedge\ldots\wedge \psi_k, v_1\wedge\ldots \wedge v_k\rangle
=\det\left(\psi_i(v_j)\right)_{1\leq i,j\leq k}$.} between $(V^*)^{\wedge k}$ and $V^{\wedge k}$, we see that $d_{ij}=c_{i+1,j+1}=\langle\varphi_{I_i},e_{I_j}\rangle$, with the notations introduced before Lemma~\ref{lemme_2}. Thus, the determinant of $D_{11}^{(h)}$ is equal to $\det\left(\varphi_{I_i}(e_{I_j})\right)_{1\leq i,j\leq h}$.

In view of the hypotheses on $A$, we have    $$e_j= \varepsilon\, e_{j-k}+\etoile\, e_{j-k+1} +\ldots + \etoile\, e_{j-1},\qquad ({\rm where\ as\ before}\ \varepsilon:=(-1)^{k-1})$$ 
for all  $j=k+1,\ldots, 2k-1$ (as in the hypothesis of Lemma~\ref{lemme_2}).  Dually we have
       $$\varphi_j= \varepsilon\, \varphi_{j-k}+\etoile\, \varphi_{j-k+1} +\ldots + \etoile\, \varphi_{j-1},$$ 
for all  $j=k+1,\ldots, 2k-1$.  
Applying Lemma~\ref{lemme_2}, we get for $1\leq i,j\leq h$ that
\begin{eqnarray*}
    \varphi_{I_i}&=&(-1)^{i+1} \varphi_{K_i}+\etoile\, \varphi_{K_{i-1}}+\ldots+\etoile\,\varphi_{K_1},\quad{\rm and}\\
     e_{I_j}&=&(-1)^{j+1} e_{K_j}+\etoile\, e_{K_{j-1}}+\ldots+\etoile\,e_{K_1}.
  \end{eqnarray*}
Using Lemma~\ref{lemme_1}, we conclude that the above determinant is equal to 
     $$\det\left((-1)^{i+j}\varphi_{K_i}(e_{K_j})\right)_{1\leq i,j\leq h},$$
which is exactly the $h\times h$-minor of the adjoint matrix of $B$, corresponding to rows and columns going from $1$ to $h$. 

To finish the argument, we apply a result  Jacobi
stating that (in the case of matrices of determinant $1$) a minor is equal to the complementary minor of the adjoint matrix.
\end{proof}

\begin{proof}[\bf Proof of Proposition~\ref{derive_equation}.]
Proposition~\ref{Condense} implies  (\ref{condensation_devive}). This equation, for $r=k$ and $s=0$ implies that the dual is a $\SL{k}$-tiling. For $r=k-1$ and $s=1$, it implies that the tiling coincides with its bidual, up to the necessary translation. 

In order to show that the dual is tame, we proceed as follows. Observe that for any matrix (finite or infinite) $(a_{ij})$ of rank at most $k$, there exist a vector space $E$ of dimension at most $k$, vectors $u_j\in E$, and linear forms $\varphi_i$ on $E$, all such that $a_{ij}=\varphi_i(u_j)$ (take the space spanned by the columns and the linear function obtained by projections of the columns). Conversely, such a data gives a matrix $(a_{ij})$ of rank at most $k$.

Now we form the matrix $(\langle \varphi_I,u_J\rangle)_{IJ}$,  over some family of $k-1$-subsets of the row and column indices. Then $u_j$ is in the $(k-1)$-th exterior power of $E$, which is of dimension at most $k$. Hence this new matrix is of rank at most $k$. This implies that the dual is at rank at most $k$, since it is obtained from the original tiling by such a construction.
\end{proof}
  Our proof of Proposition~\ref{theo_neuf} relies on the following two lemmas.
   
   \begin{lemma}\label{thm9lemma1} \bleu{Let $p$ and $q$ be two points that are adjacent horizontally, i.e.: $p=(a,b)$ and $q=(a,b+1)$.
   Then $w_q=w_p\,x\,y^i$, where $i+1$ is the number of points lying on the path $\pi_w$ in the same vertical as $q$.}
   \end{lemma}
   
 \begin{proof}[\bf Proof.]
Denote by $q_0,\ldots,q_i$ these $i+1$ points (starting from the top), and by  $r_1,\ldots,r_j$ all the points of $\pi_w$ lying to the left of $q_i$ (labelled from left to right). This corresponds to the portion of the path $\pi_w$ illustrated in Figure~\ref{portion}.
\begin{figure}[h]
$$  \setlength{\unitlength}{6mm}\setlength{\carrelength}{5mm} \begin{picture}(4,5)(0,-1.5)
   \Case{1,3}{$\!q_0$}    \Case{1.9,3}{} 
   \put(.9,1.5){$\vdots$} 
     \Case{-3,0}{$\!r_1$}       \put(-1.8,0){$\cdots$}     \Case{0.1,0}{$\!r_j$}   \Case{1,0}{$\!q_i$}  
      \Case{-3,-1}{}
\end{picture}$$
\caption{Points of the path $\pi_w$ that lie on the same row and column as $q$.}\label{portion}
\end{figure}
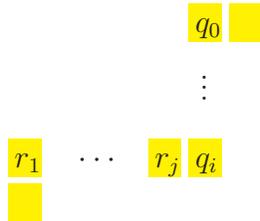
Clearly we have $w_{q_k}=y^k$, for $k<i$, and $w_{q_i}=x^j y^i$. Each of the following cases is clear (it helps to consider Figure~\ref{portion}), using the definition of $w_p$ in Section~\ref{words_tilings}.
\begin{enumerate}
\item[1)] If $q$ lies strictly above $q_0$, we have
    $$w_p=w_q\,\overline{y}^i\, \overline{x},$$
implying that $w_q=w_p\, x\, y^i$ as required. 
\item[2)]  When $q=q_j$, for $0\leq j\leq i-1$, then
    $$w_p= \overline{y}\,^{i-j} \overline{x}=y^j\overline{y}^i\overline{x}$$
 so that again we have $w_q=y^j = w_p\,x\,y^i$.
 \item[3)] Finally when $q$ lies below $q_i$,  we evidently have $w_q=w_p\,x\,y^i$, thus the assertion is verified for all possible cases.
\end{enumerate}
\end{proof}

\begin{lemma}\label{thm9lemma2}
\bleu{For $m_1,\ldots,m_k$ in $\N$, let 
     $$u_i:=x^{m_{i-1}} y\ldots x^{m_1}y, $$
 when $0\leq i\leq k$. {\rm (}In particular $u_0=1$.{\rm )}
 Then we have
    $$\begin{pmatrix} {\bf e}_k \,\mu(u_0)\\
                 {\bf e}_k \,\mu(u_1)\\
                  \vdots\\
                  {\bf e}_k\,\mu(u_{k-1})
           \end{pmatrix} = 
           \begin{pmatrix}
                0 & 0 &\ldots &0& 1\\
                0 & 0 &\ldots &1& \etoile\\
                \vdots & \vdots &\ddots & \vdots& \vdots\\
                1 & \etoile &\ldots &\etoile& \etoile\\
           \end{pmatrix}$$
}\end{lemma}
\begin{proof}[\bf Proof.]
We recursively show that ${\bf e}_k\,\mu(u_i)=(0,\ldots,0,1,\etoile,\ldots,\etoile)$ with $1$ sitting in position $(k-i)$.
If $i=0$, we have ${\bf e}_k \,\mu(u_0)={\bf e}_k=(0,\ldots,0,1)$, so that $1$ indeed sits in position $k$.
By induction we may assume that the first nonzero value of the vector
    $${\bf v}:={\bf e}_k\, \mu(x^{m_{i-1}} y\ldots x^{m_2}y)$$ 
    is a $1$ sitting in position $(k-i)$.  Then the the first nonzero value of the vector
   \begin{eqnarray*}
      {\bf e}_k\ \mu(x^{m_{i-1}} y\ldots x^{m_1})&=& {\bf v}\, \mu(x^{m_1})\\ 
                                   &=& (0,\ldots,0,1,\etoile,\ldots,\etoile),
\end{eqnarray*}
also sits in position $(k-i)$, since $\mu(x^{m_1})$ is upper unitriangular. We can thus easily conclude
since  ${\bf e}_k\,\mu(u_i)$ is obtained by multiplying (on the right) this last vector by $\mu(y)=\Id+N^{\rm tr}$,
hence its first nonzero value lies in position $(k-i-1)$.
\end{proof}

\begin{proof}[\bf Proof of Proposition~\ref{theo_neuf}.]
Let us first check that ${\mathcal T}_w$ is indeed a $\SL{k}$-tiling. Consider any set of points $p_{ij}$, $0\leq i,j\leq k-1$, forming an adjacent $k\times k$ sub-array of $\Z\times \Z$, and let us write $w_{ij}$ for the projection word $w_{p_{ij}}$ associated to these points $p_{ij}$. From Lemma~\ref{thm9lemma1} and its symmetric statement, there exists integers $m_0,\ldots ,m_{k-1}$ and $n_0,\ldots,n_{k-1}$ such that 
    $$w_{ij}= u_i\,w_{00}\,v_j,$$
with
   $$u_i= x^{m_{i-1}} y\ldots x^{m_1}y,\qquad {\rm and}\qquad v_j=x\,y^{n_1} \ldots x\,y^{m_{j-1}}.$$
We have the matrix identity
 \begin{eqnarray*}
     \Big({\bf e}_k\  \mu(u_i\,w_{00}\,v_j)\ {\bf e}_k^{\rm tr}\Big)_{0\leq i,j\leq k-1}
        &=& \Big({\bf e}_k\  \mu(u_i)\,\mu(w_{00})\, \mu(v_j)\ {\bf e}_k^{\rm tr}\Big)_{0\leq i,j\leq k-1}\\       
        &=& \begin{pmatrix} {\bf e}_k \ \mu(u_0)\\
                  \vdots\\
                  {\bf e}_k\ \mu(u_{k-1})
           \end{pmatrix}\,
            \mu(w_{00})\,
           \big( {\bf e}_k \ \mu(v_0),\ 
                  \ldots\   
                  ,\ {\bf e}_k\ \mu(v_{k-1}\big).
  \end{eqnarray*}
The fact that $\det(w_{00})=1$, together with Lemma~\ref{thm9lemma2}, implies that this matrix has determinant $1$ as announced. 

To show that $\mathcal{T}_w$ is tame we argue as follows. First, observe that any given row of $\mathcal{T}_w$ is of the form
    $$({\bf e}_k\ \mu(m)\ {\bf e}_k^{\rm tr})_{m\in \mathcal{M}},$$
where $\mathcal{M}$ is the (ordered) set of projection words of the points on this row.
Choose $k$ other rows, lying below the given row. These $k$ rows are (successively) of the form
     $$({\bf e}_k\ \mu(m_i\cdots m_1\,m)\ {\bf e}_k^{\rm tr})_{m\in \mathcal{M}},$$
 for $i$ running from $1$ to $k$, and suitable words $m_1$, $\ldots$, $m_k$. Now, the $k+1$ row vectors   ${\bf e}_k\, \mu(m_i\cdots m_1)$, $i=0,\ldots, k$, are perforce linearly dependent, since 
 they are all of length $k$. Multiplying, this linear combination by $\mu(m)\,{\bf e}_k^{\rm tr}$ on the right, we find that the $k+1$ chosen rows of $\mathcal{T}_w$ are linearly dependent, and hence  $\mathcal{T}_w$ is of rank $\leq k$.
 
 The proof that $\mathcal{T}_w$ is the same $\SL{k}$-tiling as the one described in Proposition~\ref{prop_word_tiling}, using Lemma~\ref{lem_unique}, is left to the reader.
\end{proof}

\section {Applications}\label{sec_applications}
\subsection{$\SL{2}$-Frieze patterns revisited}\label{revisit}
The aim of this section is to show that the frieze patterns of Coxeter may be realized in terms of $\SL{2}$-tilings. This gives a new slant on their study, with emphasis on their link with representations of  $\SL{2}(\Z)$.

\begin{proposition} \label{extension}
\bleu{Let $a_i$, $i\in\Z$ be non-zero elements in the field $K$. There exists a unique tame $\SL{2}$-tiling that extends the partial tiling of {\rm (\ref{frieze_esquisse})}}. 
\begin{equation}\label{frieze_esquisse}\bleu{
\begin{array}{rrrrrrrrrrrrrrrrrrrrrr} 
                        &&&\ddots\\
                         &&&\ddots& 1\\
                         &&&& a_{-1} & 1 \\
                         &&&& &a_{0} & 1 \\
                         &&&& & &a_{1} & 1 \\
                         &&&& & &&\ddots & \ddots  \\
                         &&&& & && & a_i& 1  \\
                         &&&&&&&&&\ddots&\ddots
                         \end{array}
}\end{equation}
\end{proposition}
To better study such tilings, let us consider the notion of \defn{signed continuant polynomials} $q_n(x_1,\ldots,x_n)$ defined by the recurrence
   \begin{equation}\label{continuants}
      \bleu{
         q_n(x_1,\ldots,x_n):=x_n\,  q_{n-1}(x_1,\ldots,x_{n-1})- q_{n-2}(x_1,\ldots,x_{n-2}), 
               }\end{equation}
 whenever $n>0$, setting $q_{-1}:=0$ and $q_0:=1$. We omit indices when possible, writing simply $q(x_1,\ldots,x_n)$ for $q_n(x_1,\ldots,x_n)$.                       
 Let us now consider the particular $\SL{2}$ matrices
    $$\bleu{Y(t):=\begin{pmatrix} 0 & -1\\ 1 & t \end{pmatrix}},$$
for which one may easily show by induction that 
 \begin{equation}\label{matricecontinuants}
     \bleu{Y(x_1) Y(x_2) \cdots Y(x_n)= \begin{pmatrix} -q(x_2,\ldots,x_{n-1}), & -q(x_2,\ldots,x_{n})\\
              q( x_1,\ldots,x_{n-1}), & q(x_1,\ldots,x_{n})\end{pmatrix}}.
   \end{equation}

\begin{proof} [\bf Proof  of Proposition~\ref{extension}.]
To prove unicity, we exploit the fact that the $\SL{2}$-tiling contains subarrays of the form
\begin{equation}\label{sous_matrice}
  \bleu{
\begin{matrix}
1&0&*\\
a_i&1&0
\end{matrix}}
\end{equation}
Indeed, this follows directly from the $\SL{2}$-property.
Let $C_1$, $C_2$, and $C_3$  be the three corresponding columns, from left to right. 
Then, since the tiling is of rank 2, we have $C_1-\alpha C_2+C_3=0$, which forces $\alpha=a_i$. Thus the coefficients of linearization are completely determined. Moreover, we are given at least one adjacent $2\times 2$-subarray, namely lower left $2\times 2$-submatrix of (\ref{sous_matrice}). Thus, unicity of the tiling follows by Proposition~\ref{coefflinear}. 

For the existence of the tiling, we check that we may define its entries as follows:
   $$  \bleu{\begin{matrix}   
      0 & -1 & -a_{\rouge{i}} &&\ldots & & \beta\\
      1 & 0 & -1  &\ddots &&&\vdots\\
      a_{\rouge{i}} & 1 & 0 &\ddots&\ddots&&\vdots  \\[3pt]
     \vdots & \ddots&\ddots & \ddots &\ddots&\ddots&\vdots\\
     \vdots&&\ddots&\ddots&0&-1& -a_{\rouge{j}}\\
     \vdots&&&\ddots&1 & 0 & -1\\
     \alpha & \cdots&\cdots &\cdots&a_{\rouge{j}} & 1 & 0           \end{matrix}}$$
with $\beta=-\alpha$, and  
\begin{equation}\label{valeurfrieze}
\bleu{ \alpha =q(a_i,\ldots,a_j)}.
\end{equation}
We then need only verify that the resulting tiling has the right the linearization coefficients (as in the first part of the proof). To this aim, let us denote by $\alpha'$ and $\alpha''$ the two entries of  the tiling that sit immediately to the right of $\alpha$, so that we have 
$$
\begin{matrix}
\alpha & \alpha' & \alpha''
\end{matrix}
$$
and therefore $\alpha'=q(a_{i+1},\ldots,a_j)$, and $\alpha''=q(a_{i+2},\ldots,a_j)$. But the recurrence~(\ref{continuants}) implies that
    $${\alpha-a_i\,\alpha'}+\alpha''=0.$$
Hence, since the proof for $\beta$ is analoguous, the tiling defined above has the desired linearization coefficients.
\end{proof}

Consider now any frieze pattern, as below, with the $a_i$ positive integers, having $n$ diagonals (see Figure~\ref{fig1_5}).
\begin{equation}\label{frieze}
\begin{array}{ccccccccccccccc}
&&&&&&&\ddots\\
&&&&&&&\ddots&1 \\
&&&&&&&&a_{-2}&1\\
&&&&&&&\revddots&&a_{-1}&1\\
&&&&\ddots&&\revddots&&&&a_0&1\\
&&&&&1&&&&&&a_1&1\\
&&&&&&1&&&&&&a_2&1\\
&&&&&&&1&&&&\revddots&&\ddots&\ddots\\
&&&&&&&&1&&\revddots\\
&&&&&&&&&1\\
&&&&&&&&&&\ddots\\
\end{array}
\end{equation}
Observe that, since the coefficients of the frieze pattern are positive, they are completely characterized by the $a_i$, in view of the $\SL{2}$-property. Hence this frieze pattern extends uniquely to the same complete $\SL{2}$-tiling $\A$ as the one given by Proposition~\ref{extension}. Note that this extension has a few values more, immediately deduced from the positivity of the entries of the frieze pattern, and the $\SL{2}$-property, without resulting to tameness. These are the $0$'s and $-1$'s given below. We may therefore extract from the tiling $\A$ the following subarray, where $n$ is the number of diagonals of the frieze pattern and $i\in\Z$:
$$
\begin{array}{ccccccccccccccccccc}
a_{i+1}\\
*&a_{i+2}\\
\vdots&\ddots&\ddots\\
\vdots&&&\ddots\\
1&\hdots&\hdots&*&a_{i+n-1}\\
0&1&\hdots&\hdots&*&a_{i+n}\\
-1&0&\hdots&\hdots&\hdots&*&a_{i+n+1}
\end{array}
$$
By Formula~(\ref{valeurfrieze}) for the entries of the tiling $\A$, we obtain 
$$\begin{array}{lll}
   q(a_{i+1},\ldots,a_{i+n+1})=-1, & q(a_{i+1},\ldots,a_{i+n})=0,\\
   q(a_{i+2},\ldots,a_{i+n})=1,& q(a_{i+2},\ldots,a_{i+n+1})=0
 \end{array},$$
and, using (\ref{matricecontinuants}), we conclude that   
\begin{equation}\label{=-1}
     \bleu{Y(a_{i+1}) Y(a_{i+2}) \cdots Y(a_{i+n+1})= \begin{pmatrix} -1 & 0\\
              0 & -1\end{pmatrix}}.
\end{equation}
Using this we may prove 
\begin{corollary}\label{antiperiod}
\bleu{Let  $\A=(a_{ij})$ be the unique tame $\SL{2}$-tiling $\A=(a_{ij})$ extending a given frieze pattern $\mathcal{F}$ with $n$ diagonals}.
\begin{enumerate}
\item[(i)] \bleu{$\A$ has diagonal period $n+1$, that is to say $a_{i+n+1,j+n+1}=a_{ij}$}.
\item[(ii)] {\rm (Coxeter \cite{coxeter}, Conway-Coxeter\cite{conway})} \bleu{In particular, the frieze pattern has diagonal period $n+1$}. 
\item[(iii)] \bleu{ Moreover, $\A$ has horizontal and vertical skew-period $n+1$, this is to say that}
   $$\bleu{a_{i+n+1,j}=-a_{ij}=a_{i,j+n+1}}.$$
\item[(iv)] {\rm (Coxeter \cite{coxeter}, Conway-Coxeter\cite{conway})} \bleu{ Finally, $\A$ and  $\mathcal{F}$ are invariant under  a glided symmetry, which is the symmetry with respect to the middle diagonal of $\mathcal{F}$ followed by the diagonal translation of length $\frac{n+1}{2}$}. 
\end{enumerate}
\end{corollary}


\begin{proof}[\bf Proof.]
Referring to (37), let $C_i$ denote the column containing the coefficient $a_i$. Then we have $C_{i}-a_iC_{i+1}+C_{i+2}=0$, as is shown at the beginning of the proof of Proposition~\ref{extension}. Thus we have the following recurrence between the $\Z\times 2$ matrices $(C_i,C_{i+1})$:
$$
(C_i,C_{i+1})=(C_{i+1},C_{i+2})\, Y(a_i).
$$
Thus, (\ref{=-1}) implies that $\A$ has horizontal skew-period $n+1$. Vertical periodicity follows by symmetry, and the diagonal periodicity follows at once.
In order to prove (iv), note that  (\ref{=-1}) implies 
$$
Y(a_{i+2}) \cdots Y(a_{i+n+1})=-Y(a_{i+1}^{-1})= \begin{pmatrix} -a_{i+1} & -1\\
              1 & 0\end{pmatrix}.
$$
Thus $a_{i+1}=q(a_{i+3},\ldots,a_{i+n})$ by  (\ref{matricecontinuants}). This shows, by taking $i=n,0,1,\ldots$ and recalling that we have the diagonal period $n+1$ (hence $a_i=a_{i+n+1}$) that: $a_{n+1}=q(a_{n+3},\ldots, a_{2n})=q(a_2,\ldots, a_{n-1})$, $a_1=q(a_3,\ldots,a_n)$, \ldots.
Hence, using (\ref{valeurfrieze}) 
we see
 that  $R$ has the following form, extending  (\ref{frieze_esquisse}):
$$
\bleu{
\begin{array}{cccccccccccccccccccccc} 
                        &&&\ddots\\
                         &&&\ddots& 1\\
                         &&&& a_{1} & 1 \\
                         &&&& &a_{2} & 1 \\
                         &&&& & &a_{3} & 1 \\
                         &&&\ddots&\ddots & &&\ddots & \ddots  \\
                         &&&& 1& a_{n+1}&& & a_{n-1}& 1  \\
                         &&&&&1&a_1&&&a_n&1\\
                         &&&&&&1&a_2&&&a_{n+1}&1\\
                         &&&&&&&\ddots&\ddots&&&\ddots&\ddots
                         \end{array}
}
$$
We conclude by using a symmetric version of Proposition~\ref{extension}.
\end{proof}

Following Conway-Coxeter (in \cite{conway}) we call \defn{quiddity} a sequence $a_1,\ldots,a_{n+1}$, where $a_i$ gives the number of triangle incident to the vertex $i$ in a triangulation of a convex $(n+1)$-gone, whose vertex are successively labeled $1$ to $n+1$ turning around the $n$-gone. They show \cite[p. 180]{conway} that any quiddity may be obtained from the particular quiddity $1\,1\,1$ by successive applications of the local rewriting rule
$$   \bleu{\cdots\,\rouge{a\,b}\,\cdots \quad\rightarrow\quad \cdots\,\rouge{a+1\,1\,b+1}\,\cdots}.$$

We prove below their result that quiddities and frieze patterns are in one-to-one correpondence. For this, we make a detour through presentations of the group $\SL{2}(\Z)$. 

\begin{proposition}\label{rewriting}
\bleu{
Consider the rewriting rule in the free monoid $\mathbb{P}^*$ generated by $\mathbb P$}
\begin{equation}
  \bleu{ (a+1)\,\,1\,\,(b+1) \quad \rightarrow \quad a\,b},
\end{equation}
\bleu{where $a,b \in \mathbb P$. Then 
   $$Y(w):=Y(n_1)\ldots Y(n_k)=\varepsilon\,\begin{pmatrix} 1 & 0 \\ 0 & 1\end{pmatrix},\qquad (\varepsilon =\pm 1),$$
if and only if $w\rightarrow^* 1^k$, with $k\equiv 0\ (\mod 6)$ when $\varepsilon =1$, and $k\equiv 3\ (\mod 6)$ when $\varepsilon =-1$. In this case, if $w$ is not a power of 1, then it contains a factor $(a+1)\,\,1\,\,(b+1)$.
}
\end{proposition}

One direction of the proposition easily follows from the identities
\begin{equation}\label{relations1}
   \bleu{Y(1)^3=-\begin{pmatrix} 1 & 0 \\ 0 & 1\end{pmatrix}},\qquad {\rm and}\qquad \bleu{Y(a+1)Y(1)Y(b+1)=Y(a)Y(b)},
\end{equation}
both of which can be easily checked by direct computation. 

Moreover, for further use, it is easily checked that
\begin{equation}\label{relation2}
\bleu{Y(1)Y(2)Y(1)Y(2)=-1}
\end{equation}
and also, recursively, that 
    \begin{equation}\label{Yn}
       \bleu{Y(n)=(-1)^n (Y(2)\,Y(1)^2)^{n-2}\,Y(2)}.
    \end{equation}
    
We now give a proof of Proposition~\ref{rewriting} after recalling some facts regarding presentations of $\SL{2}(\Z)$. To simplify our discussion, let us informally\footnote{This can easily be made formally correct by adding a generator, with straightforward relations, to our presentations.}   write ``$-1$'' for a central element of $\SL{2}$ whose square is the identity (denoted by $1$).

\begin{lemma}
\bleu{Denoting $Y(1)$ by $y_1$, and $Y(2)$ by $y_2$}, 
\begin{enumerate}
\item[(i)] \bleu{$\SL{2}(\Z)$ affords the presentation }
\begin{equation}\label{first_pres}
   \bleu{ \langle\ y_1, y_2\ |\ y_1^3=-1,\quad (y_1y_2)^2=-1\ \rangle.}
\end{equation}
\item[(ii)] \bleu{$\SL{2}(\Z)$ affords the confluent presentation}
\begin{equation}\label{conf_pres}
   \bleu{ \langle\ y_1, y_2\ |\  y_1^3\rightarrow -1,\quad y_2y_1y_2\rightarrow y_1^2\ \rangle.}
 \end{equation}
 \end{enumerate}
\end{lemma}

\begin{proof}[\bf Proof.]
To show (i), let $$
a=\begin{pmatrix}1&1\\0&1 \end{pmatrix}, b=\begin{pmatrix}1&0\\1&1 \end{pmatrix}.
$$
It is well-known that $\SL{2}(\Z)$ has the presentation 
\begin{equation}\label{presentationSL2}
\bleu{ \langle\ a,b\ |\ \bar a b\bar a=b\bar a b,(\bar a b\bar a)^4=1\rangle}.
\end{equation}
Direct calculations show that 
\begin{equation}\label{ab_en_fonction}
a=\bar y_1y_2,\quad b=y_1 \bar y_2.
\end{equation}
Thus $y_1,y_2$ generate $\SL{2}(\Z)$. The relations in (i) hold by  (\ref{relations1}) and (\ref{relation2}). Therefore, it is enough to show that these relations imply the relations in  (\ref{presentationSL2}), once $a,b$ have been replaced using  (\ref{ab_en_fonction}). 
By direct substitution, we get 
    $$\bar a b\bar a=\bar y_2 y_1y_1\bar y_2\bar y_2y_1,\quad{\rm and}\quad b\bar a b=y_1\bar y_2\bar y_2y_1y_1\bar y_2.$$
     Now, since $y_1^3=-1$, we have $-\bar y_1=y_1^2$, and hence $y_1^2\bar y_2^3(-\bar y_1)=(-\bar y_1)\bar y_2^3y_1^2$. Thus $y_1^2\bar y_2^2(-\bar y_2\bar y_1)=(-\bar y_1\bar y_2)\bar y_2^2y_1^2$. But, since $y_1y_2y_1y_2=-1$, we also have  $y_2y_1y_2y_1=-1$, and therefore $-\bar y_2\bar y_1=y_1y_2$ and $-\bar y_1\bar y_2=y_2y_1$. Hence 
$$ y_1^2\bar y_2^2y_1y_2=y_2y_1\bar y_2^2y_1^2.$$
Multiplying this both on the left and on the right by $\bar y_2$, we obtain
$$\bar y_2y_1^2\bar y_2^2y_1=y_1\bar y_2^2y_1^2\bar y_2,$$
so that $\bar a b\bar a=b\bar a b$.

On the other hand, we have $\bar a b=\bar y_2y_1^2\bar y_2$, and we have seen that $y_1y_2=-\bar y_2\bar y_1$. Thus
$$y_1=-\bar y_2\bar y_1\bar y_2=-\bar y_2(-y_1^2)\bar y_2=\bar y_2y_1^2 \bar y_2,$$
since $\bar y_1=-y_1^2$. Thus $(\bar y_2y_1^2\bar y_2)^6=1$. It follows, using  $\bar a b\bar a=b\bar a b$, that
$$ (\bar a b\bar a)^4=\bar a b\bar a b\bar a b \bar a b\bar a b\bar a b=(\bar a b)^6=(\bar y_2y_1^2\bar y_2)^6=1.$$

(ii) We conclude from the first part that
\begin{equation}
   \bleu{ \langle\ y_1, y_2\ |\  y_1^3= -1,\quad y_2y_1y_2=y_1^2\ \rangle.}
 \end{equation}
is a presentation of $\SL{2}(\Z)$. Orienting the equalities, we obtain a rewriting system, whose confluence we must prove. This follows since
 the only non-trivial critical pair that needs to be examined is
$$\begin{matrix}   & & y_2y_1y_2y_1y_2\\
       &  \swarrow & & \searrow\\
     (y_1^2)\,y_1y_2   && &&  y_2y_1\, (y_1^2)\\
     &\searrow && \swarrow\\
         && -y_2
      \end{matrix}$$
   This ends our proof.   
\end{proof}



\begin{proof}[\bf Proof of Proposition~\ref{rewriting}.] 
We need only show that if $w\not=1^n$ is such that $Y(w)=\pm Id$, then $w$ must contain a factor of the form $(a+1)\,\,1\,\,(b+1)$. To see this, formally replace each letter $n$ in ${\mathbb P}$ by $y_n$ in words $w$ in ${\mathbb P}^*$. Then, assuming that $w$ is different from $y_1^n$, we may consider the canonical expansion 
       $$w= y_1^{n_0}\,y_{m_1}\,y_1^{n_1}\,y_{m_2}\,\cdots y_{m_k}\,y_1^{n_k},$$
where each $m_i\geq 2$ and $k\geq 1$. Using (\ref{Yn}), we replace in this expansion each $y_m$ by $(-1)^m\, (y_2y_1^2)^{m-2} y_2$, we obtain a word in $\{y_1,y_2\}^*$ containing at least one instance of $y_2$. Since the system~(\ref{conf_pres}) is confluent, this word must contain $y_2y_1y_2$, hence one of the $n_i$ must be equal to $1$, thus proving our assertion.
\end{proof}

\begin{corollary}[Conway-Coxeter \cite{conway}]
\bleu{For each frieze pattern of the form {\rm (}\ref{frieze}{\rm )}, the bi-infinite sequence of positive integers $\cdots a_{-2}a_{-1}a_0a_1a_2a_3\cdots$ is equal to $\cdots wwwwwww \cdots$ for some quiddity $w$}.
\end{corollary}

\begin{proof}[\bf Proof.]
Denote by $\mathcal{F}$ this frieze pattern, let $n$ be the number of diagonals of $\mathcal{F}$ and denote by $\A$ the $\SL{2}$-tiling obtained through Proposition~\ref{extension}. Then, by the discussion before Corollary~\ref{antiperiod}, the coefficients $a_i$ satisfy 
$$\bleu{Y(a_{1}) Y(a_{2}) \cdots Y(a_{n+1})= \begin{pmatrix} -1 & 0\\
              0 & -1\end{pmatrix}}
.$$
Thus, by Corollary~\ref{antiperiod}, $\cdots a_{-2}a_{-1}a_0a_1a_2a_3\cdots=\cdots wwwwwww \cdots$, with $w=a_1\cdots a_{n+1}$. Moreover, by Proposition~\ref{rewriting}, we have 
     $$\bleu{a_{i-1}>1,\quad  a_i=1},\quad {\rm and}\quad  \bleu{a_{i+1}>1},$$
 for some $i$, $1<i<n+1$. Thus we find in $\A$  the subarray
$$
\begin{array}{ccccccccc}
a_{i-2}&1&0\\
&a_{i-1}&1&0\\
&(a_{i-1}-1)&1&1&0\\
&&(a_{i+1}-1)&a_{i+1}&1\\
&&&&a_{i+2}
\end{array}
$$
Let $C_j$ and $R_j$ denote the columnn and row containing $a_j$. Then, as in the discussion at the beginning of the proof of Proposition~\ref{rewriting}, we have $C_i-C_{i+1}+C_{i+2}=0$ and $R_i-R_{i-1}+R_{i-2}=0$. Thus, by Lemma~\ref{suppression_colonne}, we may suppress both the column $C_{i+1}$ and  the row $R_{i-1}$, to get the tame $\SL{2}$-tiling
$$
\begin{array}{cccccccc}
a_{i-2}&1\\
&(a_{i-1}-1)&1&\\
&&(a_{i+1}-1)&1\\
&&&a_{i+2}
\end{array}
$$
We may clearly do this periodically for each column $C_{i+1+p(n+1)}$ and each row $R_{i-1+p(n+1)}$, for $p\in \Z$. Since both $\A$ and $\mathcal{F}$ have the diagonal period $n+1$, by Corollary~\ref{antiperiod}, we obtain a tame $\SL{2}$-tiling $\A'$ and a frieze pattern $\mathcal{F}'$ with $n-1$ diagonals, such that $\A'$ is the unique extension of $\mathcal{F}'$ according to Proposition~\ref{extension}. This proves the corollary, once it is noted that the initial case corresponds to the frieze patterns reduced to $n=2$ diagonals containing only $1$'s (here considered as having a diagonal period equal to $3$).
\end{proof}

\subsection{$\N\times \N$ $\SL{k}$-tilings}\label{NbyN}
When we restrict ourselves to $\N\times \N$ arrays, we may apply tools from matrix algebra and generating series.
Assume that, for an invertible $k\times k$ matrix,  we have a $(k+h)\times (k+h)$ matrix of rank $k$  that decomposes into blocks in the following manner   
\begin{equation}\begin{pmatrix}\label{matrice_bloc}
            \mathcal{S} & \Lambda\\
            \Gamma & \mathcal{X}
            \end{pmatrix},
\end{equation}
with $h$ possibly infinite.
Then we must necessarily have $\bleu{\mathcal{X}=\Gamma\, \mathcal{S}^{-1}\,\Lambda}$.
Indeed, it is clear in the following  simple matrix identity     
$$\begin{pmatrix}
            \mathcal{S} & \Lambda\\
            \Gamma & \mathcal{X}
            \end{pmatrix}\,
            \begin{pmatrix}
            \Id_k & -\mathcal{S}^{-1}\,\Lambda\\
            0 & Id_h
            \end{pmatrix}=
            \begin{pmatrix}
            \mathcal{S} & 0\\
            \Gamma & \mathcal{X}-\Gamma\, \mathcal{S}^{-1}\,\Lambda
            \end{pmatrix}
        $$
that the right-hand side is also a matrix of rank $k$, since we are multiplying our original rank $k$ matrix by an invertible one. However, we already know that $\mathcal{S}$ is of rank $k$. This forces $\mathcal{X}-\Gamma\, \mathcal{S}^{-1}\,\Lambda$ to vanish, and we have the required identity. 

Let us assume that $\A$  is a tame \defn{quarter-plane} (of shape $\N\times \N$) $\SL{k}$-array. Choose $\mathcal{S}:=\A_{00}^{(k)}$, and let $\Gamma$ (resp. $\Lambda$) stand for the subarray consisting of the first $k$ columns (resp. rows) of $\A$. Then, we deduce from the above identity that  we have
\begin{equation}\label{decomposition}
  \bleu{   \A=\Gamma\, \mathcal{S}^{-1}\,\Lambda},
\end{equation} 
whenever $\A$ and $\mathcal{S}$ are both of rank $k$. 
It follows that for any subset $I$ (resp. $J$) of rows (resp. columns, with $\#I=\#J$), we have 
\begin{equation}
    \bleu{  \A_{IJ}=\Gamma_I\, \mathcal{S}^{-1}\, \Lambda_J}.
  \end{equation}
  The simplest possible case of this identity allows the calculation of entries of $\A$ in the form
  \begin{equation}\label{entrees_de_A}
       \bleu{a_{ij}= \Gamma_i\, \mathcal{S}^{-1}\, \Lambda_j}.
   \end{equation}
 Now, if we choose both $I$ and $J$ to be of cardinality $k$, and take the determinant of both sides, we deduce from the fact that $\det(\mathcal S)=1$, the
identity
\begin{equation}
    \bleu{\M_{IJ}=\M_{I,\{1,\ldots,k\}} \, \M_{\{1,\ldots,k\},J}}.
  \end{equation}


A straightforward encoding of the $\N\times \N$-array $\A$ is through its bivariate \defn{generating function}:
  \begin{equation}\label{gen_fonct}
         \bleu{\A(x,y):=\sum_{(i,j)}  a_{ij}\, x^i y^j},
   \end{equation}
with the sum running over all pairs $(i,j)$ belonging to the shape of $\A$. An equivalent description may be given in terms of matrices, considering $X=\left(x^i\right)_{0\leq i}$  as an infinite one-line matrix, and likewise
$Y=\left(y^j\right)_{0\leq j}$ as an infinite one-column matrix. We then have
    $\A(x,y)=X\,\A\,Y$.
Now, when $\A$ is tame, it follows from (\ref{decomposition}) that 
we have
\begin{eqnarray}
    \A(x,y)&=&X\,\A\,Y \nonumber\\
             &=& X\, \Gamma\, \mathcal{S}^{-1}\, \Lambda\, Y\nonumber\\
             &=& \bleu{\begin{pmatrix} C_1(x) & \ldots & C_k(x) \end{pmatrix}
                    \, \mathcal{S}^{-1}\, \begin{pmatrix} L_1(y) \\ \vdots \\ L_k(y)\end{pmatrix}},\label{gen_fonct_matrice}
    \end{eqnarray}
 where the $C_i(x)$ are respectively the generating functions of the first $k$ columns of $\A$. Likewise
 the $L_j(x)$ are the respective generating functions of the first $k$ rows of $\A$. Often, as below, we have $\mathcal{S}={\rm Id}$.

To illustrate the situation considered above, one may show that the $\SL{k}$-property holds for the matrix of binomial coefficients 
     $$\Gamma:=\left(\binom{j}{i}\right)_{\monatop{0\leq i,}{0\leq j <k}}.$$ 
From this, we get a $\SL{k}$-array $\A:=\Gamma\Lambda$, with $\Lambda$ equal to the transpose of $\Gamma$. Observe that
the generating function of the $j^{\rm th}$-column (resp. $i^{\rm th}$) of $\Gamma$ (resp. $\Lambda$) is evidently
      $$\frac{1}{(1-x)^j} = \sum_{i\geq 0} \binom{i+j}{i} x^i,\qquad {\rm for}\ j=0,1,\ldots, k-1$$
(resp. $(1-y)^i$). After calculation, using (\ref{gen_fonct_matrice}), we get that the generating function of $\A$ is
 \begin{equation}\label{gen_fonct_k}
     \bleu{\A(x,y)=\sum_{\ell=1}^{k} \frac{x^{\ell-1} y^{\ell-1}}{(1-x)^\ell (1-y)^\ell}}.
  \end{equation}
The following result follows from the constructions in Section~\ref{words_tilings}.
\begin{proposition}\label{prop_rectangle}
\bleu{The tiling given by {\rm (}\ref{gen_fonct_k}{\rm )} has all minors of the from $\M_{i0}^{(m)}$ and $\M_{0j}^{(m)}$ equal to $1$, whenever $m< k$, and it is a $\SL{k}$-tiling.}
\end{proposition}

Using (\ref{entrees_de_A}), or directly from (\ref{gen_fonct_k}), one may calculate that the individual entries of $\A$ are given by the formula
 \begin{equation}\label{coeff_rect}
    \bleu{a_{ij}=\sum_{\ell=0}^{k-1} \binom{i}{\ell}\binom{j}{\ell}}.
 \end{equation}
It follows also from Section~\ref{words_tilings}, that for $(i,j)$ in the $k$-fringe, $a_{ij}=\binom{i+j}{i}$ .
For example, with $k=3$, we get the array of Figure~\ref{fig_quarter}.
\begin{figure}[ht]
$$
\begin{picture}(0,0)(0,0)\setlength{\unitlength}{5mm}
\put(1.7,3.5){\put(0,0){\vector(1,0){17}}
             \put(0,0){\vector(0,-1){6.5}}}
\end{picture}
\A=   \begin {array}{cccccccccc}
1&1&1&1&1&1&1&1&1&\cdots\\ 
1&2&3&4&5&6&7&8&9&\cdots\\ 
1&3&6&10&15&21&28&36&45&\cdots\\ 
1&4&10&19&31&46&64&85&109&\cdots\\ 
1&5&15&31&53&81&115&155&201&\cdots\\
\vdots & \vdots &\vdots &\vdots &\vdots &\vdots &\vdots &\vdots &\vdots &\ddots
\end {array}
$$
\caption{A $\SL{3}$-tiling of $\N\times \N$.} \label{fig_quarter}
\end{figure} 

It may readily be shown that the dual of $\A$ affords the generating function
 \begin{equation}\label{gen_fonct_derk}
     \bleu{\A^*(x,y)=\frac{1}{(1-x)(1-y)}+\sum_{\ell=2}^k \frac{x\,y}{(1-x)^\ell (1-y)^\ell}}.
  \end{equation}

\subsection{Zigzag path}
\begin{figure}[ht]
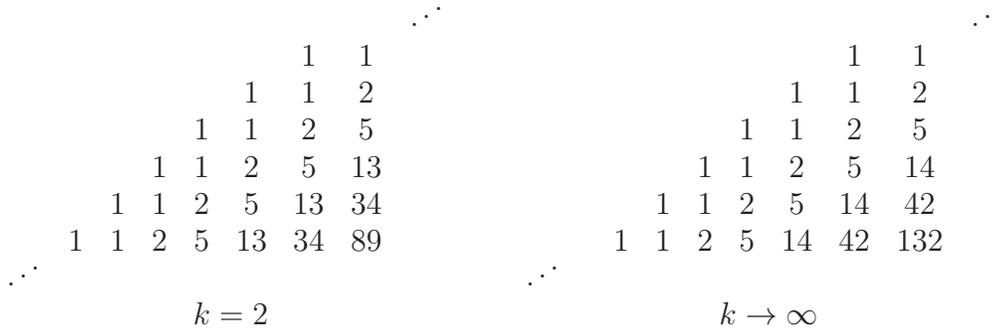

$$
\begin{array}{c}
\begin{array}{ccccccccccccccccccccccccccccccccccccccccccccccc}
&&&&&&&&\revddots&&&&&&&&&&&&\revddots\\
&&&&&&1&1&&&&&&&&&&&1&1\\
&&&&&1&1&2&&&&&&&&&&1&1&2\\
&&&&1&1&2&5&&&&&&&&&1&1&2&5\\
&&&1&1&2&5&13&&&&&&&&1&1&2&5&14\\
&&1&1&2&5&13&34&&&&&&&1&1&2&5&14&42\\
&1&1&2&5&13&34&89&&&&&&1&1&2&5&14&42&132\\
\revddots&&&&&&&&&&&\revddots
\end{array}\\
k=2 \hskip6cm k\rightarrow \infty
\end{array}
$$\vskip-10pt
\caption{Zigzag path tiling}\label{fig_zigzag}
 \end{figure}

It is shown in \cite{assem} that the $\SL{2}$-tiling associated to the bi-infinite word $\cdots xyxyxyxyx \cdots$ has entries equal to the Fibonacci numbers of even rank (if we set $F_{n+2}=F_{n+1}+F_n$, $F_0=F_1=1$), see the left part of Figure~\ref{fig_zigzag}. If we let $k$ go to infinity, then by Section~\ref{words_tilings}, the entries of the resulting tiling are the Catalan numbers. In particular, it is noteworthy that the value of the $(k\times k)$-principal minors given by Proposition~\ref{prop_word_tiling}  corresponds in this situation to the classical result stating that for any natural integer $k$, the Hankel matrix $(C_{h+i+j})_{i,j=0,\ldots, k}$ (with either $h=0$, or $h=1$) has determinant equal to $1$ (here, as usual, we have $C_n=\frac{1}{n+1}\binom{2n}{n}$), see the right part of Figure~\ref{fig_zigzag}.

\section{Closing remarks}\label{sec_closing}

\subsection*{A converse}
Experiments suggest that a ``converse'' of Proposition~\ref{derive_equation} holds, namely that for any tame tiling, if {\rm (}\ref{condensation_devive}{\rm )} holds for some pair $(r,s)$ for which $r+s=k$, then the tiling is necessarily a $\SL{k}$-tiling. Special cases, for small values of $k$, are easy to prove using generic value tilings and Gr\"obner basis computations.

\subsection*{Generalized frieze patterns}
A notion of generalized frieze patterns, for $k>2$, has been considered in \cite{cordes}.  These are best understood in terms of certain tame ``toric''  $\SL{k}$-tilings $\mathcal{A}$. 
More precisely, we say that a tiling has a \defn{skew-period} $(p,q)$ in $\Z\times \Z$ ($\not=(0,0)$), if and only if 
   \begin{equation}\label{defn_skew_periodic}
        \mathcal{A}(i+p,j+q) =(-1)^k\mathcal{A}(i,j),\qquad {\rm for\ all}\ (i,j)\in\Z\times \Z,
     \end{equation}
and we then say that the tiling is \defn{skew-periodic}. A \defn{toric} tiling is one that has two linearly independent skew-periods. A \defn{$\SL{k}$-frieze patterns} $\mathcal{A}$ is a tame $\SL{k}$-tiling such that
\begin{equation}\label{defn_cond_frise}
   \mathcal{A}(i,j)= \begin{cases}
      1, & \text{if}\  i=j,\\
      0, & \text{if}\ i-j<k,\\
      (-1)^{k-1}, & \text{if}\ i=j+k,
\end{cases}
\end{equation}
which is periodic (not skew) with a period of the form $(p,-p)$, for $p>k$.
In other words, on top of being periodic, the tiling is prescribed to have a diagonal of $1$'s, another diagonal filled with $(-1)^{k-1}$ below, with these two diagonals separated by  $(k-1)$ diagonals of $0$'s.

Condition (\ref{defn_cond_frise}) and periodicity (together with tameness) ensure that the whole tiling is determined by its values along a band $\{(i,j)\  j\leq i< j+p,\ j\in\Z\}$, with $p$ as above. Moreover, one may show that any such tiling is toric, with skew-periods $(p,0)$ and $(0,p)$. This implies that it exhibits a frieze-like behaviour, since the tiling must necessarily have period $(p,p)$. The generalized frieze patterns of \cite{cordes} appear as special cases of this notion. 

It is interesting to observe (and easy to prove) that the $\SL{k}$-frieze patterns, with $p=k+2$, that are  of the form illustrated in Figure~\ref{fig_slk_frises}, include those for which the sequence $\ldots,a_{-1},a_0,a_1,\ldots a_i\,\ldots$ is a quiddity. Indeed, applying the $\SL{k}$-condition to the submatrix with main diagonal corresponding to $a_i$'s, one easily check that we have
    $$q_k(a_i,a_{i+1},\ldots,a_{i+k-1})=1,$$
 with $q_k$ denoting the signed continuant polynomials considered in Section~\ref{revisit}. 
\begin{figure}[ht]
\begin{equation}\bleu{
\begin{array}{ccccccccccccccccccc} 
                       \ddots &\ddots&\ddots&\ddots&\ddots\\
                           &0&0&1& a_{-1} & 1 \\
                         &&0&0&1&a_{0} & 1 \\
                         &&&0& 0&1&a_{1} & 1 \\
                         &&&&\ddots& \ddots&\ddots&\ddots & \ddots  \\
                         &&&& &0&0& 1& a_i& 1  \\
                         &&&&&&\ddots&\ddots&\ddots&\ddots&\ddots
                         \end{array}
}\end{equation}
\caption{Positive integer $\SL{k}$-frize patterns of ``width'' $1$.}\label{fig_slk_frises}
\end{figure}

Other interesting toric $\SL{k}$-tilings seem to abound. For example, with $k=3$ and $p=4$, we have the following positive integer valued toric $\SL{k}$-tilings:
$$\mathcal{A}=\begin {array}{cccccccccccccc} 
\ddots&\vdots&\vdots&\vdots&\vdots&\vdots&\vdots&\vdots&\vdots&\vdots&\vdots&\vdots&\revddots\\ \noalign{\medskip}
\cdots&1&1&1&1&1&1&1&1&1&1&1&\cdots\\ \noalign{\medskip}
\cdots&1&2&3&2&1&2&3&2&1&2&3&\cdots\\ \noalign{\medskip}
\cdots&1&3&6&4&1&3&6&4&1&3&6&\cdots\\ \noalign{\medskip}
\cdots&1&2&4&3&1&2&4&3&1&2&4&\cdots\\ \noalign{\medskip}
\cdots&1&1&1&1&1&1&1&1&1&1&1&\cdots\\ \noalign{\medskip}
\cdots&1&2&3&2&1&2&3&2&1&2&3&\cdots\\ \noalign{\medskip}
\cdots&1&3&6&4&1&3&6&4&1&3&6&\cdots\\ \noalign{\medskip}
\revddots&\vdots&\vdots&\vdots&\vdots&\vdots&\vdots&\vdots&\vdots&\vdots&\vdots&\vdots&\ddots\end {array}
\begin{picture}(0,0)(0,0)
\rouge{\put(-195,-4){\framebox(60,72){}}}
\end{picture}
$$
It may be checked that this is a tame tiling, 
however the entries of the corresponding dual tiling are not all positive, since:
$$\partial\mathcal{A}=\begin {array}{rrrrrrrrrrrrrrrrrrrrrrrr} 
\ddots&\vdots&\vdots&\vdots&\vdots&\vdots&\vdots&\vdots&\vdots&\vdots&\vdots&\revddots\\ \noalign{\medskip}
\cdots&1&1&-1&-1&1&1&-1&-1&1&1&\cdots\\ \noalign{\medskip}
\cdots&1&3&0&-2&1&3&0&-2&1&3&\cdots\\ \noalign{\medskip}
\cdots&-1&0&2&1&-1&0&2&1&-1&0&\cdots\\ \noalign{\medskip}
\cdots&-1&-2&1&2&-1&-2&1&2&-1&-2&\cdots\\ \noalign{\medskip}
\cdots&1&1&-1&-1&1&1&-1&-1&1&1&\cdots\\ \noalign{\medskip}
\cdots&1&3&0&-2&1&3&0&-2&1&3&\cdots\\ \noalign{\medskip}
\revddots&\vdots&\vdots&\vdots&\vdots&\vdots&\vdots&\vdots&\vdots&\vdots&\vdots&\ddots
\end {array}$$

General properties of tame toric tilings, as well as results concerning $\SL{k}$-frieze patterns similar to those of Subsection~\ref{revisit}, will be the subject of a planed sequel to this paper.

\subsection*{T-systems}
On a closing note, it is interesting to observe that there is a close tie between tame $\SL{k}$-tilings and the notion of $T$-systems, which appear as solutions of the discrete Hirota equation (see \cite{di_Francesco}) of mathematical physics. Indeed, up to a simple relabelling, one can characterize the entries of $T$-systems in terms of derivatives of suitably chosen tame $\SL{k}$-tilings. Recall that a   $T$-system $T:\{0,\ldots, r+1\}\times\Z\times \Z\rightarrow \N$ must satisfy the equation
   \begin{equation}\label{T_def}
         T_{\{\alpha,j,k+1\}}T_{\{\alpha,j,k-1\}}=T_{\{\alpha,j+1,k\}}T_{\{\alpha,j-1,k\}}+T_{\{\alpha+1,j,k+1\}}T_{\{\alpha-1,j,k+1\}},
    \end{equation}
with boundary conditions
     \begin{equation}\label{T_cond}
         T_{\{0,j,k\}}=T_{\{r+1,j,k\}}=1,
    \end{equation}
 for all $j$, and $k$ in $\Z$. It is shown in \cite{di_Francesco} that
       \begin{equation}\label{T_formule}
                T_{\{\alpha,j,k\}}= \det \begin{pmatrix}  
                                                      T_{\{1,j-a+b,k+a+b-\alpha-1\}}
                                                \end{pmatrix}  _{1\leq a,b\leq \alpha}.
     \end{equation}
From this, one can readily see that 
     $$T_{\alpha,j,k}=(\partial^\alpha \mathcal{A})_{s,t}$$
 for $s$, $t$ simple linear functions of $j$ and $k$, and $\mathcal{A}$ a $\SL{r+1}$-tiling directly obtained from $(T_{1,j,k})_{j,k}$.



\begin{thebibliography}{10}   

\bibitem{assem} 
\auteur{I.~Assem, C.~Reutenauer, and D.~Smith},
\titreref{Frises}, To appear.

\bibitem{baker}
\auteur{Andrew Baker},
\titreref{Matrix Groups, An Introduction to Lie Group Theory},
Springer-Verlag, 2002.

\bibitem{bressoud}
\auteur{D.M. Bressoud},
\titreref{Proofs and Confirmations, The Story of the Alternating Matrix Conjecture},
Cambridge University Press, 1999.

\bibitem{conway}
\auteur{J.~Conway and H.S.M.~Coxeter},
\titreref{Triangulated polygons and frieze patterns},
The Mathematical Gazette \vol{57} (1973), 87--94 and 175-183.

\bibitem{cordes}
\auteur{C.M.~Cordes and D.P.~Roselle},
\titreref{Generalized Frieze Patterns},
Duke Math. J. \vol{39}, Number 4 (1972), 637--648.

\bibitem{coxeter}
\auteur{H.S.M.~Coxeter},
\titreref{Frieze Patterns}, Acta Arithmetica XVIII (1971), 297--310.

\bibitem{di_Francesco} 
\auteur{Ph. Di Francesco, and R. Kedem},
\titreref{Positivity of the $T$-system cluster algebra},
arXiv:0908.3122.


\bibitem{alice}
\auteur{Rev. C. L. Dodgson},
\titreref{Condensation of Determinants, being a brief Method for computing their arithmetical values},
Proceedings of the Royal Society XV (1866),  150--155.


\bibitem{fomin} 
\auteur{S.~Fomin and A.~Zelevinsky},
\titreref{The Laurent Phenomenon},
Advances in Applied Mathematics,
\vol{28} (2002), 119--144.  

\bibitem{viennot}
\auteur{I.~Gessel and G.X.~Viennot}



\bibitem{propp} 
\auteur{J.~Propp},
\titreref{The combinatorics of frieze patterns and Markoff numbers},\\
Proceedings of FPSAC'06, also available at
arXiv:math/0511633v4 [math.CO].


\end{thebibliography}
\end{document}